\def\draft{n}

\documentclass[12pt]{amsart}
\usepackage{amssymb,epic,eepic,xr}
\usepackage[dvips, final]{graphicx}
\usepackage{fullpage}
\theoremstyle{plain}

\newtheorem{theorem}{Theorem}
\newtheorem{proposition}{Proposition}[section]

\newtheorem{lemma}[proposition]{Lemma}

\theoremstyle{definition}
\newtheorem{definition}[proposition]{Definition}

\theoremstyle{remark}

\newtheorem{remark}[proposition]{Remark}

\def\printname#1{
        \if\draft y
                \smash{\makebox[0pt]{\hspace{-0.5in}
                        \raisebox{8pt}{\tt\tiny #1}}}
        \fi
}

\newcommand{\mathmode}[1]{$#1$}

\newlength{\standardunitlength}
\catcode`\@=11
\font\bfcaptionfont=cmssbx10 scaled \magstephalf
\newdimen\@captionmargin\@captionmargin=2\parindent
\long\def\@makecaption#1#2{
    \vskip 10pt
    \setbox\@tempboxa\hbox{%
      \small\sf{\bfcaptionfont #1. }\ignorespaces #2}
    \ifdim \wd\@tempboxa >\captionwidth {
        \rightskip=\@captionmargin\leftskip=\@captionmargin
        \unhbox\@tempboxa\par}
      \else
        \hbox to\hsize{\hfil\box\@tempboxa\hfil}
    \fi}
\catcode`\@=12

\def\lbl#1{\label{#1}\printname{#1}}

\def\qed{{\hfill\text{$\Box$}}}

\newenvironment{myitemize}{
        \begin{list}{$\bullet$}{\setlength{\leftmargin}{16pt}
        \setlength{\labelwidth}{12pt}
        \setlength{\labelsep}{4pt}}
}{
        \end{list}
}

\newlength{\globalparindent}
\def\mathcenter#1{%
  \vcenter{\hbox{#1}}%
}

\def\graph#1{
        \includegraphics[trim=-2 -2 -2 -2]{#1}
}

\def\mathgraph#1{
        \mathcenter{\graph{#1}}
}

\def\mathgraphsmall#1{
        \mathcenter{\includegraphics{#1}}
}

\def\mathgraphint#1{
        \mathcenter{\includegraphics[trim=-1 -1 -1 -1]{#1}}
}

\newcommand{\Aa}{\textsl\AA}
\newcommand{\Ao}{{\mathcal A^\circ}}
\newcommand{\Arhus}{\AA rhus}

\newcommand{\B}{\mathcal B}

\newcommand{\Gint}{\int^{FG}}
\newcommand{\LMO}{\text{\rm LMO}}

\newcommand{\NDint}{\int^{(m)}}

\newcommand{\QQ}{\mathbb Q}
\newcommand{\RR}{\mathbb R}

\newcommand{\bbZ}{\mathbb Z}
\newcommand{\Zcheck}{\check{Z}}
\newcommand{\calA}{\mathcal A}
\newcommand{\calB}{\mathcal B}

\newcommand{\frakg}{\mathfrak g}

\newcommand{\sgn}{\mathop{\text{\rm sgn}}\nolimits}

\renewcommand{\qed}{~\hfill$\square$}

\newcommand\downstrut[2]{{\mathop\downarrow^{#1}_{#2}}}
\newcommand\strutn[2]{\,^{#1}\!\!\frown^{#2}}
\newcommand\strutu[2]{\,_{#1}\!\!\smile_{#2}}
\newcommand\strutv[2]{{\mathop\mid^{#1}_{#2}}}
\newcommand\upstrut[2]{{\mathop\uparrow^{#2}_{#1}}}

\newcommand{\doubletheta}{\mathgraphint{draws/comp.40}}
\newcommand{\crossedtet}{\mathgraphint{draws/comp.41}}

\newcommand\prtl[1]{{\partial_{#1}}}

\newcommand\comma{\mathbin ,}

\begin{document}
\newdimen\captionwidth\captionwidth=\hsize

\title[The \AA{}rhus integral III: The relation with the LMO invariant]
{
  The \AA{}rhus integral of rational homology 3-spheres III: The Relation
  with the Le-Murakami-Ohtsuki Invariant
}

\author[Bar-Natan]{Dror~Bar-Natan}
\address{
  Institute of Mathematics\\
  The Hebrew University\\
  Giv'at-Ram, Jerusalem 91904\\
  Israel
}
\curraddr{
  Department of Mathematics\\
  University of Toronto\\
  Toronto, Ontario, Canada M5S 3G3
}
\email{drorbn@math.toronto.edu}
\urladdr{http://www.math.toronto.edu/\~{}drorbn}

\author[Garoufalidis]{Stavros~Garoufalidis}
\address{Department of Mathematics\\
        Harvard University\\
        Cambridge MA 02138\\
        USA}
\curraddr{School of Mathematics\\
  Georgia Institute of Technology\\
  Atlanta, GA  30332-0160\\
  USA
}
\email{stavros@math.gatech.edu}
\urladdr{http://www.math.gatech.edu/\~{}stavros}

\author[Rozansky]{Lev~Rozansky}
\address{Department of Mathematics, Statistics, and Computer Science\\
  University of Illinois at Chicago\\
  Chicago IL 60607-7045\\
  USA}
\curraddr{Department of Mathematics, UNC-CH\\
  CB 3250 Phillips Hall\\
  Chapel Hill, NC 27599-3250\\
  USA
}
\email{rozansky@email.unc.edu}

\author[Thurston]{Dylan~P.~Thurston}
\address{Department of Mathematics\\
  University of California at Berkeley\\
  Berkeley CA 94720-3840\\
  USA}
\curraddr{Department of Mathematics\\
  Harvard University\\
  Cambridge MA 02138\\
  USA
}
\email{dpt@math.harvard.edu}

\thanks{This preprint is available electronically at
  {\tt http://www.math.toronto.edu/\~{}drorbn}, at \newline
  {\tt http://www.math.gatech.edu/\~{}stavros}, and at {\tt
    arXiv:math/9808013}.
}

\date{This edition: Sep.~15,~1903; \ \ First edition: August 4, 1998.}

\begin{abstract}
Continuing the work started in~\cite{Aarhus:I} and~\cite{Aarhus:II},
we prove the relationship between the \Arhus{} integral and the
invariant $\Omega$ (henceforth called $\LMO$) defined by T.Q.T.~Le,
J.~Murakami and T.~Ohtsuki in~\cite{LeMurakamiOhtsuki:Universal}. The
basic reason for the relationship is that both constructions afford
an interpretation as ``integrated holonomies''. In the case of the
\Arhus{} integral, this interpretation was the basis for everything we did
in~\cite{Aarhus:I} and~\cite{Aarhus:II}. The main tool we used there was
``formal Gaussian integration''. For the case of the $\LMO$ invariant,
we develop an interpretation of a key ingredient, the map $j_m$, as
``formal negative-dimensional integration''. The relation between the
two constructions is then an immediate corollary of the relationship
between the two integration theories.
\end{abstract}

\maketitle

\tableofcontents

\section{Introduction}

This paper is the third in a four-part series on ``the \Arhus{}
integral of rational homology 3-spheres''. In Part~I of this
series~\cite{Aarhus:I}, we gave the definition of a diagram-valued
invariant $\Aa$ of rational homology spheres.\footnote{
  A precise definition of $\Aa$ appears in~\cite{Aarhus:I}. It
  is a good idea to have~\cite{Aarhus:I} as well
  as~\cite{LeMurakamiOhtsuki:Universal} handy while reading this paper,
  as many of the definitions introduced and explained in those articles
  will only be repeated here in a very brief manner.
}
In Part~II (\cite{Aarhus:II}) we proved that $\Aa$ is a well-defined
invariant of rational homology 3-spheres and that it is universal in the
class of finite type invariants of integral homology spheres. In this
paper we show that $\Aa$, when defined, is essentially equal
to the invariant $\LMO$ defined earlier by Le, Murakami, and Ohtsuki
in~\cite{LeMurakamiOhtsuki:Universal}.

Both invariants $\LMO$ and $\Aa$ take values in the space
$\calA(\emptyset)$, the completed graded space of manifold diagrams
modulo the $AS$ and $IHX$ relations, as defined in detail
in~\cite[Definition~\ref{I-def:ManifoldDiagrams}]{Aarhus:I} and as
recalled briefly in Figure~\ref{fig:ManifoldDiagrams}, and the precise
statement of their near-equality is as follows:

\begin{figure}
\[
  \printname{VacuumCC}
        \setlength{\unitlength}{0.5\standardunitlength}
        \begin{array}{c}  \hspace{-1.7mm}
                \raisebox{-8pt}{\begingroup\makeatletter\ifx\SetFigFont\undefined
\def\x#1#2#3#4#5#6#7\relax{\def\x{#1#2#3#4#5#6}}%
\expandafter\x\fmtname xxxxxx\relax \def\y{splain}%
\ifx\x\y   
\gdef\SetFigFont#1#2#3{%
  \ifnum #1<17\tiny\else \ifnum #1<20\small\else
  \ifnum #1<24\normalsize\else \ifnum #1<29\large\else
  \ifnum #1<34\Large\else \ifnum #1<41\LARGE\else
     \huge\fi\fi\fi\fi\fi\fi
  \csname #3\endcsname}%
\else
\gdef\SetFigFont#1#2#3{\begingroup
  \count@#1\relax \ifnum 25<\count@\count@25\fi
  \def\x{\endgroup\@setsize\SetFigFont{#2pt}}%
  \expandafter\x
    \csname \romannumeral\the\count@ pt\expandafter\endcsname
    \csname @\romannumeral\the\count@ pt\endcsname
  \csname #3\endcsname}%
\fi
\fi\endgroup
{\renewcommand{\dashlinestretch}{30}
\begin{picture}(2420,1256)(0,-10)
\put(608,612){\ellipse{1200}{1200}}
\put(2108,612){\ellipse{600}{900}}
\path(158,987)(608,612)(1058,987)
\path(608,612)(608,12)
\path(1808,612)(2408,612)
\path(1433,987)(1883,912)
\path(1883,312)(1133,312)
\path(908,1137)	(961.376,1167.463)
	(1008.459,1191.462)
	(1050.126,1209.215)
	(1087.258,1220.944)
	(1151.428,1227.204)
	(1208.000,1212.000)

\path(1208,1212)	(1265.377,1178.903)
	(1317.052,1135.747)
	(1361.696,1084.256)
	(1397.982,1026.150)
	(1424.584,963.150)
	(1440.172,896.977)
	(1443.420,829.354)
	(1433.000,762.000)

\path(1433,762)	(1410.875,714.085)
	(1369.516,675.457)
	(1303.649,642.601)
	(1259.877,627.238)
	(1208.000,612.000)

\end{picture}
} }
                \hspace{-1.9mm}
        \end{array}

  \qquad\qquad\printname{IHX}
        \setlength{\unitlength}{0.5\standardunitlength}
        \begin{array}{c}  \hspace{-1.7mm}
                \raisebox{-8pt}{\begingroup\makeatletter\ifx\SetFigFont\undefined
\def\x#1#2#3#4#5#6#7\relax{\def\x{#1#2#3#4#5#6}}%
\expandafter\x\fmtname xxxxxx\relax \def\y{splain}%
\ifx\x\y   
\gdef\SetFigFont#1#2#3{%
  \ifnum #1<17\tiny\else \ifnum #1<20\small\else
  \ifnum #1<24\normalsize\else \ifnum #1<29\large\else
  \ifnum #1<34\Large\else \ifnum #1<41\LARGE\else
     \huge\fi\fi\fi\fi\fi\fi
  \csname #3\endcsname}%
\else
\gdef\SetFigFont#1#2#3{\begingroup
  \count@#1\relax \ifnum 25<\count@\count@25\fi
  \def\x{\endgroup\@setsize\SetFigFont{#2pt}}%
  \expandafter\x
    \csname \romannumeral\the\count@ pt\expandafter\endcsname
    \csname @\romannumeral\the\count@ pt\endcsname
  \csname #3\endcsname}%
\fi
\fi\endgroup
{\renewcommand{\dashlinestretch}{30}
\begin{picture}(2724,639)(0,-10)
\put(1738,235){\makebox(0,0)[lb]{\smash{{\mathmode{-}}}}}
\path(12,612)(612,612)
\path(312,612)(312,12)
\path(12,12)(612,12)
\path(1062,612)(1062,12)
\path(1062,312)(1662,312)
\path(1662,612)(1662,12)
\path(2112,612)(2712,12)
\path(2712,612)(2112,12)
\path(2262,162)(2562,162)
\put(687,237){\makebox(0,0)[lb]{\smash{{\mathmode{=}}}}}
\end{picture}
} }
                \hspace{-1.9mm}
        \end{array}

  \qquad\qquad\printname{AS}
        \setlength{\unitlength}{0.5\standardunitlength}
        \begin{array}{c}  \hspace{-1.7mm}
                \raisebox{-8pt}{\begingroup\makeatletter\ifx\SetFigFont\undefined
\def\x#1#2#3#4#5#6#7\relax{\def\x{#1#2#3#4#5#6}}%
\expandafter\x\fmtname xxxxxx\relax \def\y{splain}%
\ifx\x\y   
\gdef\SetFigFont#1#2#3{%
  \ifnum #1<17\tiny\else \ifnum #1<20\small\else
  \ifnum #1<24\normalsize\else \ifnum #1<29\large\else
  \ifnum #1<34\Large\else \ifnum #1<41\LARGE\else
     \huge\fi\fi\fi\fi\fi\fi
  \csname #3\endcsname}%
\else
\gdef\SetFigFont#1#2#3{\begingroup
  \count@#1\relax \ifnum 25<\count@\count@25\fi
  \def\x{\endgroup\@setsize\SetFigFont{#2pt}}%
  \expandafter\x
    \csname \romannumeral\the\count@ pt\expandafter\endcsname
    \csname @\romannumeral\the\count@ pt\endcsname
  \csname #3\endcsname}%
\fi
\fi\endgroup
{\renewcommand{\dashlinestretch}{30}
\begin{picture}(2029,639)(0,-10)
\path(12,612)(312,312)(612,612)
\path(312,312)(312,12)
\path(1362,312)(1362,12)
\path(1062,612)	(1110.175,601.033)
	(1154.410,590.038)
	(1194.844,578.933)
	(1231.613,567.634)
	(1294.707,544.128)
	(1344.791,518.861)
	(1410.324,460.408)
	(1437.000,387.000)

\path(1437,387)	(1416.615,336.139)
	(1362.000,312.000)

\path(1362,312)	(1307.385,336.139)
	(1287.000,387.000)

\path(1287,387)	(1313.670,460.408)
	(1379.201,518.861)
	(1429.286,544.128)
	(1492.381,567.634)
	(1529.152,578.933)
	(1569.587,590.038)
	(1613.824,601.033)
	(1662.000,612.000)

\put(687,237){\makebox(0,0)[lb]{\smash{{\mathmode{+}}}}}
\put(1737,237){\makebox(0,0)[lb]{\smash{{\mathmode{=\ 0}}}}}
\end{picture}
} }
                \hspace{-1.9mm}
        \end{array}

\]
\caption{
  A connected manifold diagram of degree $6$ (half the number of
  vertices in it) and the $IHX$ and $AS$ relations.
} \lbl{fig:ManifoldDiagrams}
\end{figure}
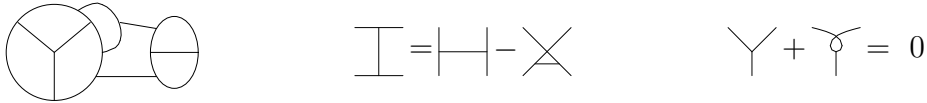

\begin{theorem} \lbl{thm:main} (Proof in Section~\ref{sec:dull})
Let $M$ be a rational homology sphere and let $|H^1(M)|$ denote the
number of elements in its first cohomology group (over $\bbZ$).  Then
\begin{equation} \lbl{eq:main}
  \Aa(M) = |H^1(M)|^{-\deg} \LMO(M),
\end{equation}
where in the graded space $\calA(\emptyset)$, $|H^1(M)|^{-\deg}$
denotes the operation that multiplies any degree $m$ element by
$|H^1(M)|^{-m}$.
\end{theorem}

In particular, if $M$ is an integral homology sphere, i.e., if
$|H^1(M)|=1$, then simply $\Aa(M) = \LMO(M)$. Also note that
Equation~\eqref{eq:main} implies that $\Aa(M)=\hat\Omega(M)$,
where $\hat\Omega$ is the invariant defined
in~\cite[Section~6.2]{LeMurakamiOhtsuki:Universal}.

The definitions of $\LMO$ and of $\Aa$ are very similar. Let us trace
this similarity to the point where the two definitions diverge. That
point is of course the key point of our paper, for it is there that we
have to prove something non-trivial. In reading the following few
paragraphs, on the similarity and differences between the definitions of
$\LMO$ and of $\Aa$, the reader may find it helpful to consult with
Figure~\ref{fig:diagram} which summarizes the maps and the spaces involved.

\begin{figure}
\[ {
  \def\R{{\left\{\parbox{0.4in}{\begin{center}
    \scriptsize regular pure tangles
  \end{center}}\right\}}}
  \def\Z{{\scriptstyle\Zcheck}}
  \def\LMMO{{\parbox{1.1in}{\begin{center}\scriptsize
    the LMMO \\ version of the \\ Kontsevich integral
  \end{center}}}}
  \def\AX{{\calA(\uparrow_X)}}
  \def\s{{\scriptstyle\sigma}}
  \def\f{{\parbox{0.5in}{\begin{center}
    \scriptsize formal PBW
  \end{center}}}}
  \def\B{{\calB(X)}}
  \def\FG{{\scriptstyle\Gint}}
  \def\j{{\scriptstyle j_m}}
  \def\ndi{{\scriptstyle\NDint}}
  \def\Aar{{\scriptstyle\Aa_0}}
  \def\Am{{\calA(\emptyset)}}
  \def\L{{\scriptstyle\LMO_0^{(m)}}}
  \def\Ar{{\Ao(\emptyset) / (O_m, P_{m+1})}}
  \def\isom{{\parbox{0.6in}{\begin{center}
    \scriptsize isomorphic in degrees $\leq m$
  \end{center}}}}
  \printname{diagram}
        \setlength{\unitlength}{0.75\standardunitlength}
        \begin{array}{c}  \hspace{-1.7mm}
                \raisebox{-8pt}{\begingroup\makeatletter\ifx\SetFigFont\undefined%
\gdef\SetFigFont#1#2#3#4#5{%
  \reset@font\fontsize{#1}{#2pt}%
  \fontfamily{#3}\fontseries{#4}\fontshape{#5}%
  \selectfont}%
\fi\endgroup%
{\renewcommand{\dashlinestretch}{30}
\begin{picture}(7615,2400)(0,-10)
\put(5812.500,1237.500){\arc{2476.136}{5.6700}{6.8964}}
\path(6863.615,642.511)(6825.000,525.000)(6914.371,610.515)
\path(6914.371,1864.485)(6825.000,1950.000)(6863.615,1832.489)
\path(5100,1050)(6300,450)
\path(6179.252,476.833)(6300.000,450.000)(6206.085,530.498)
\path(3675,1200)(4650,1200)
\path(4530.000,1170.000)(4650.000,1200.000)(4530.000,1230.000)
\path(1350,1200)(2850,1200)
\path(2730.000,1170.000)(2850.000,1200.000)(2730.000,1230.000)
\path(5100,1425)(6300,1950)
\path(6202.086,1874.417)(6300.000,1950.000)(6178.037,1929.386)
\path(675,1725)(676,1725)(677,1726)
	(681,1728)(686,1731)(693,1735)
	(703,1740)(715,1746)(731,1754)
	(749,1763)(770,1774)(794,1786)
	(821,1799)(851,1814)(883,1829)
	(918,1844)(954,1861)(993,1878)
	(1033,1894)(1074,1911)(1117,1928)
	(1161,1944)(1206,1960)(1252,1975)
	(1299,1990)(1348,2004)(1398,2017)
	(1449,2029)(1502,2041)(1557,2051)
	(1615,2061)(1675,2070)(1738,2079)
	(1805,2086)(1875,2094)(1950,2100)
	(1994,2103)(2040,2107)(2088,2110)
	(2137,2113)(2189,2116)(2243,2119)
	(2299,2121)(2358,2124)(2419,2126)
	(2483,2129)(2550,2131)(2620,2133)
	(2692,2135)(2767,2137)(2845,2139)
	(2926,2141)(3009,2143)(3095,2144)
	(3184,2146)(3276,2147)(3369,2149)
	(3465,2150)(3564,2152)(3664,2153)
	(3766,2155)(3869,2156)(3974,2157)
	(4080,2158)(4186,2159)(4293,2160)
	(4401,2161)(4508,2162)(4614,2163)
	(4719,2164)(4823,2165)(4926,2166)
	(5026,2167)(5124,2168)(5219,2168)
	(5310,2169)(5399,2170)(5483,2170)
	(5564,2171)(5640,2171)(5711,2172)
	(5778,2172)(5840,2173)(5897,2173)
	(5949,2173)(5996,2174)(6038,2174)
	(6075,2174)(6107,2174)(6135,2174)
	(6158,2175)(6177,2175)(6192,2175)
	(6204,2175)(6213,2175)(6225,2175)
\path(6105.000,2145.000)(6225.000,2175.000)(6105.000,2205.000)
\path(675,675)(676,675)(677,674)
	(681,672)(686,669)(693,665)
	(703,660)(715,654)(731,646)
	(749,637)(770,626)(794,614)
	(821,601)(851,586)(883,571)
	(918,556)(954,539)(993,522)
	(1033,506)(1074,489)(1117,472)
	(1161,456)(1206,440)(1252,425)
	(1299,410)(1348,396)(1398,383)
	(1449,371)(1502,359)(1557,349)
	(1615,339)(1675,330)(1738,321)
	(1805,314)(1875,306)(1950,300)
	(1994,297)(2040,293)(2088,290)
	(2137,287)(2189,284)(2243,281)
	(2299,279)(2358,276)(2419,274)
	(2483,271)(2550,269)(2620,267)
	(2692,265)(2767,263)(2845,261)
	(2926,259)(3009,257)(3095,256)
	(3184,254)(3276,253)(3369,251)
	(3465,250)(3564,248)(3664,247)
	(3766,245)(3869,244)(3974,243)
	(4080,242)(4186,241)(4293,240)
	(4401,239)(4508,237)(4614,237)
	(4719,236)(4823,235)(4926,234)
	(5026,233)(5124,232)(5219,232)
	(5310,231)(5399,230)(5483,230)
	(5564,229)(5640,229)(5711,228)
	(5778,228)(5840,227)(5897,227)
	(5949,227)(5996,226)(6038,226)
	(6075,226)(6107,226)(6135,226)
	(6158,225)(6177,225)(6192,225)
	(6204,225)(6213,225)(6225,225)
\path(6105.000,195.000)(6225.000,225.000)(6105.000,255.000)
\path(3300,975)(3301,974)(3303,972)
	(3306,968)(3311,962)(3318,953)
	(3328,942)(3340,928)(3354,912)
	(3371,893)(3390,873)(3411,852)
	(3434,830)(3458,808)(3484,785)
	(3512,763)(3540,742)(3571,722)
	(3602,702)(3636,684)(3672,668)
	(3711,652)(3752,638)(3797,624)
	(3846,612)(3900,600)(3931,594)
	(3964,588)(3999,582)(4036,576)
	(4076,570)(4117,564)(4162,558)
	(4209,552)(4259,546)(4311,540)
	(4366,534)(4424,528)(4485,522)
	(4549,515)(4615,509)(4683,503)
	(4753,496)(4826,490)(4900,483)
	(4976,476)(5053,470)(5130,463)
	(5208,457)(5286,450)(5364,444)
	(5440,437)(5515,431)(5587,425)
	(5658,420)(5725,414)(5789,409)
	(5849,404)(5905,400)(5957,396)
	(6004,392)(6046,389)(6083,386)
	(6115,383)(6143,381)(6166,380)
	(6184,378)(6199,377)(6209,376)(6225,375)
\path(6103.362,352.544)(6225.000,375.000)(6107.105,412.427)
\put(0,1125){\makebox(0,0)[lb]{\smash{{\mathmode{\R}}}}}
\put(3750,900){\makebox(0,0)[lb]{\smash{{\mathmode{\f}}}}}
\put(1200,750){\makebox(0,0)[lb]{\smash{{\mathmode{\LMMO}}}}}
\put(4725,1125){\makebox(0,0)[lb]{\smash{{\mathmode{\B}}}}}
\put(6300,225){\makebox(0,0)[lb]{\smash{{\mathmode{\Ar}}}}}
\put(6300,2025){\makebox(0,0)[lb]{\smash{{\mathmode{\Am}}}}}
\put(3000,0){\makebox(0,0)[lb]{\smash{{\mathmode{\L}}}}}
\put(2925,1125){\makebox(0,0)[lb]{\smash{{\mathmode{\AX}}}}}
\put(5775,1500){\makebox(0,0)[lb]{\smash{{\mathmode{\FG}}}}}
\put(4050,1275){\makebox(0,0)[lb]{\smash{{\mathmode{\s}}}}}
\put(4800,600){\makebox(0,0)[lb]{\smash{{\mathmode{\j}}}}}
\put(2025,1275){\makebox(0,0)[lb]{\smash{{\mathmode{\Z}}}}}
\put(3075,2250){\makebox(0,0)[lb]{\smash{{\mathmode{\Aar}}}}}
\put(5775,825){\makebox(0,0)[lb]{\smash{{\mathmode{\ndi}}}}}
\put(7200,1200){\makebox(0,0)[lb]{\smash{{\mathmode{\isom}}}}}
\end{picture}
} }
                \hspace{-1.9mm}
        \end{array}

} \]
\caption{
  The main ingredients in the definitions of $\LMO_0$ and of $\Aa_0$:
  The top half of this diagram is the definition of $\Aa_0$. The bottom
  half shows the two definitions of $\LMO_0$ --- the original involving
  $j_m$ and our variant using $\NDint$; the two are equivalent by
  Lemma~\ref{lem:equivalence}. Finally, the triangle on the right is
  nearly commutative, as detailed in Proposition~\ref{prop:main}.
} \lbl{fig:diagram}
\end{figure}
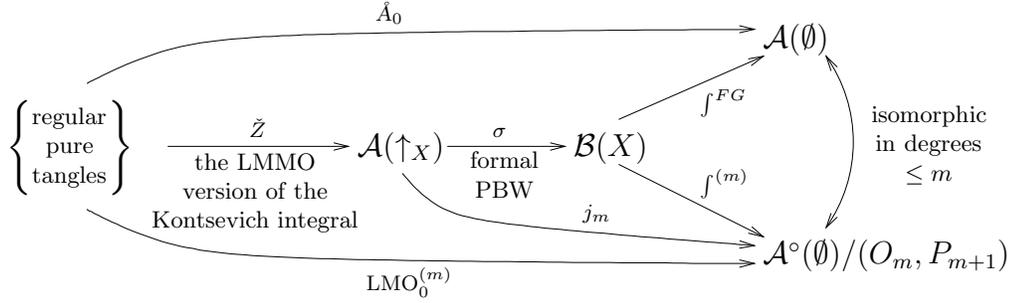

The definitions of $\LMO$ and of $\Aa$ start with the definitions of
their ``pre-normalized'' versions $\LMO_0$ and $\Aa_0$, which are
invariants of regular links (framed links having a non-degenerate
linking matrix) that are also invariant under the second Kirby move and
both use the same renormalization procedure to ensure invariance also
under the first Kirby move. The resulting $\LMO$ and $\Aa$ are therefore
invariants of rational homology spheres presented by surgery over regular
links.

The definitions of $\LMO_0$ and of $\Aa_0$ are also similar.  For both
it is beneficial to start with regular pure tangles (framed pure
tangles having a non-degenerate linking matrix,
see~\cite[Definition~\ref{I-def:RPT}]{Aarhus:I}) rather than with
regular links. By closure, every regular pure tangle defines a regular
link and every regular link is obtained in this way. Both $\LMO_0$ and
$\Aa_0$ descend from invariants of regular pure tangles to invariants
of links; $\LMO_0$ is defined that way
in~\cite{LeMurakamiOhtsuki:Universal} to start with, and for $\Aa_0$ it
is shown in~\cite[Section~\ref{II-subsec:DescendsToLinks}]{Aarhus:II}.

Both $\LMO_0$ and $\Aa_0$ are defined as compositions of several
maps. In both cases the first map is $\Zcheck$, the Kontsevich
integral in its
Le-Murakami-Murakami-Ohtsuki~\cite{LeMurakami2Ohtsuki:3Manifold}
normalization (check~\cite[Definition~\ref{I-def:Zcheck}]{Aarhus:I} for
the adaptation to pure tangles). If $X$ is the set of components of a
given regular pure tangle (or more elegantly, a set of ``labels'' or
``colors'' for these components), the first map $\Zcheck$ takes its
values in $\calA(\uparrow_X)$, the completed graded space of chord
diagrams for $X$-labeled pure tangles (modulo the usual $4T/STU$
relations; see a precise definition
in~\cite[Definition~\ref{I-def:TangleDiagrams}]{Aarhus:I} and a brief
reminder in Figure~\ref{fig:Adef}).

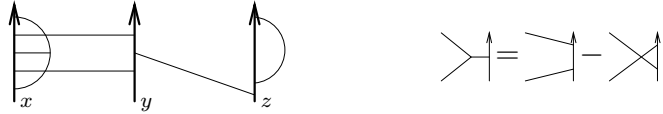
\begin{figure}
\[
  \printname{PTDiagram}
        \setlength{\unitlength}{0.5\standardunitlength}
        \begin{array}{c}  \hspace{-1.7mm}
                \raisebox{-8pt}{\begingroup\makeatletter\ifx\SetFigFont\undefined
\def\x#1#2#3#4#5#6#7\relax{\def\x{#1#2#3#4#5#6}}%
\expandafter\x\fmtname xxxxxx\relax \def\y{splain}%
\ifx\x\y   
\gdef\SetFigFont#1#2#3{%
  \ifnum #1<17\tiny\else \ifnum #1<20\small\else
  \ifnum #1<24\normalsize\else \ifnum #1<29\large\else
  \ifnum #1<34\Large\else \ifnum #1<41\LARGE\else
     \huge\fi\fi\fi\fi\fi\fi
  \csname #3\endcsname}%
\else
\gdef\SetFigFont#1#2#3{\begingroup
  \count@#1\relax \ifnum 25<\count@\count@25\fi
  \def\x{\endgroup\@setsize\SetFigFont{#2pt}}%
  \expandafter\x
    \csname \romannumeral\the\count@ pt\expandafter\endcsname
    \csname @\romannumeral\the\count@ pt\endcsname
  \csname #3\endcsname}%
\fi
\fi\endgroup
{\renewcommand{\dashlinestretch}{30}
\begin{picture}(3407,1352)(0,-10)
\put(22.000,715.000){\arc{900.000}{4.7124}{7.8540}}
\put(2984.500,752.500){\arc{828.402}{4.8030}{7.7633}}
\thicklines
\path(22,115)(22,1315)
\path(82.000,1075.000)(22.000,1315.000)(-38.000,1075.000)
\path(1522,115)(1522,1315)
\path(1582.000,1075.000)(1522.000,1315.000)(1462.000,1075.000)
\path(3022,115)(3022,1315)
\path(3082.000,1075.000)(3022.000,1315.000)(2962.000,1075.000)
\thinlines
\path(472,715)(22,715)
\path(22,940)(1522,940)
\path(22,490)(1522,490)
\path(1522,715)(3022,190)
\put(97,40){\makebox(0,0)[lb]{\smash{{\mathmode{\scriptstyle x}}}}}
\put(1597,40){\makebox(0,0)[lb]{\smash{{\mathmode{\scriptstyle y}}}}}
\put(3097,40){\makebox(0,0)[lb]{\smash{{\mathmode{\scriptstyle z}}}}}
\end{picture}
} }
                \hspace{-1.9mm}
        \end{array}

  \qquad\qquad\printname{STU}
        \setlength{\unitlength}{0.5\standardunitlength}
        \begin{array}{c}  \hspace{-1.7mm}
                \raisebox{-8pt}{\begingroup\makeatletter\ifx\SetFigFont\undefined
\def\x#1#2#3#4#5#6#7\relax{\def\x{#1#2#3#4#5#6}}%
\expandafter\x\fmtname xxxxxx\relax \def\y{splain}%
\ifx\x\y   
\gdef\SetFigFont#1#2#3{%
  \ifnum #1<17\tiny\else \ifnum #1<20\small\else
  \ifnum #1<24\normalsize\else \ifnum #1<29\large\else
  \ifnum #1<34\Large\else \ifnum #1<41\LARGE\else
     \huge\fi\fi\fi\fi\fi\fi
  \csname #3\endcsname}%
\else
\gdef\SetFigFont#1#2#3{\begingroup
  \count@#1\relax \ifnum 25<\count@\count@25\fi
  \def\x{\endgroup\@setsize\SetFigFont{#2pt}}%
  \expandafter\x
    \csname \romannumeral\the\count@ pt\expandafter\endcsname
    \csname @\romannumeral\the\count@ pt\endcsname
  \csname #3\endcsname}%
\fi
\fi\endgroup
{\renewcommand{\dashlinestretch}{30}
\begin{picture}(2724,639)(0,-10)
\put(1738,236){\makebox(0,0)[lb]{\smash{{\mathmode{-}}}}}
\path(612,12)(612,612)
\path(642.000,492.000)(612.000,612.000)(582.000,492.000)
\path(1662,12)(1662,612)
\path(1692.000,492.000)(1662.000,612.000)(1632.000,492.000)
\path(2712,12)(2712,612)
\path(2742.000,492.000)(2712.000,612.000)(2682.000,492.000)
\path(12,612)(387,312)(12,12)
\path(387,312)(612,312)
\path(1062,612)(1662,462)
\path(1062,12)(1662,162)
\path(2112,612)(2712,162)
\path(2112,12)(2712,462)
\put(687,237){\makebox(0,0)[lb]{\smash{{\mathmode{=}}}}}
\end{picture}
} }
                \hspace{-1.9mm}
        \end{array}

\]
\caption{
  A diagram in $\calA(\uparrow_{\{x,y,z\}})$ and the $STU$ relation.
} \lbl{fig:Adef}
\end{figure}

Here the two definitions diverge, though a part of this divergence is
rather minor.

\begin{myitemize}

\item $\Aa_0$ is defined to be the composition
$\Aa_0=\Gint\circ\,\sigma\circ\Zcheck$. Here $\sigma$ denotes the
diagrammatic version of the Poincare-Birkhoff-Witt theorem (defined as
in~\cite{Bar-Natan:Vassiliev, Bar-Natan:Homotopy}, though normalized
slightly differently, as
in~\cite[Definition~\ref{I-def:sigma}]{Aarhus:I}) with values in
$\calB(X)$, the completed graded space of $X$-marked uni-trivalent
diagrams modulo $AS$ and $IHX$ relations (``Chinese characters''
in~\cite{Bar-Natan:Vassiliev, Bar-Natan:Homotopy}; see a precise
definition
in~\cite[Definition~\ref{I-def:UnitrivalentDiagrams}]{Aarhus:I} and a
brief reminder in Figure~\ref{fig:TypicalCC}). In~\cite{Aarhus:I,
Aarhus:II} we have discussed extensively how the space $\calB(X)$ can
be viewed as a space of functions, and how the partially defined map
$\Gint$ can be viewed as ``formal Gaussian integration'',
(see~\cite[Definition~\ref{I-def:Gint}]{Aarhus:I} and
Section~\ref{sec:proof} of this article for a definition, and
Appendix~\ref{sec:comp-gint} for an example of how to apply it). In a sense
detailed
in~\cite{Aarhus:I}, $\Gint$ is a diagrammatic analogue of the usual
notion of perturbed Gaussian integration --- it is defined by breaking
the ``integrand'' into ``quadratic'' and ``higher order'' terms,
inverting the quadratic, and gluing the higher order terms to each
other using the inverse quadratic as glue, in the spirit of Feynman
diagrams.

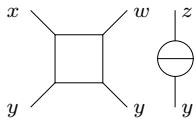
\begin{figure}
\[ \printname{TypicalCC}
        \setlength{\unitlength}{0.5\standardunitlength}
        \begin{array}{c}  \hspace{-1.7mm}
                \raisebox{-8pt}{\begingroup\makeatletter\ifx\SetFigFont\undefined
\def\x#1#2#3#4#5#6#7\relax{\def\x{#1#2#3#4#5#6}}%
\expandafter\x\fmtname xxxxxx\relax \def\y{splain}%
\ifx\x\y   
\gdef\SetFigFont#1#2#3{%
  \ifnum #1<17\tiny\else \ifnum #1<20\small\else
  \ifnum #1<24\normalsize\else \ifnum #1<29\large\else
  \ifnum #1<34\Large\else \ifnum #1<41\LARGE\else
     \huge\fi\fi\fi\fi\fi\fi
  \csname #3\endcsname}%
\else
\gdef\SetFigFont#1#2#3{\begingroup
  \count@#1\relax \ifnum 25<\count@\count@25\fi
  \def\x{\endgroup\@setsize\SetFigFont{#2pt}}%
  \expandafter\x
    \csname \romannumeral\the\count@ pt\expandafter\endcsname
    \csname @\romannumeral\the\count@ pt\endcsname
  \csname #3\endcsname}%
\fi
\fi\endgroup
{\renewcommand{\dashlinestretch}{30}
\begin{picture}(2359,1390)(0,-10)
\put(2100,715){\ellipse{450}{450}}
\path(600,1015)(1200,1015)(1200,415)
	(600,415)(600,1015)
\path(300,1315)(600,1015)
\path(1200,1015)(1500,1315)
\path(1200,415)(1500,115)
\path(600,415)(300,115)
\path(1875,715)(2325,715)
\path(2100,940)(2100,1315)
\path(2100,490)(2100,115)
\put(0,1240){\makebox(0,0)[lb]{\smash{{\mathmode{\scriptstyle x}}}}}
\put(0,40){\makebox(0,0)[lb]{\smash{{\mathmode{\scriptstyle y}}}}}
\put(1575,40){\makebox(0,0)[lb]{\smash{{\mathmode{\scriptstyle y}}}}}
\put(1575,1240){\makebox(0,0)[lb]{\smash{{\mathmode{\scriptstyle w}}}}}
\put(2175,1240){\makebox(0,0)[lb]{\smash{{\mathmode{\scriptstyle z}}}}}
\put(2175,40){\makebox(0,0)[lb]{\smash{{\mathmode{\scriptstyle y}}}}}
\end{picture}
} }
                \hspace{-1.9mm}
        \end{array}
 \]
\caption{An $\{x,y,z,w\}$-marked uni-trivalent diagram.}
\lbl{fig:TypicalCC}
\end{figure}

\item $\LMO_0$ is defined to be an ``assembly'' of maps $\LMO_0^{(m)}$,
where $m\geq 0$ is an integer. That is, the degree $m$ piece of
$\LMO_0$ is defined to be the degree $m$ piece of $\LMO_0^{(m)}$ for
each $m$. Each map $\LMO_0^{(m)}$ is a composition $j_m\circ\Zcheck$
where each $j_m$, and hence each map $\LMO_0^{(m)}$, has a different
target space. These target spaces are the spaces $\Ao(\emptyset) /
(O_m, P_{m+1})$, where $\Ao(\emptyset)$ is the same as the space
$\calA(\emptyset)$ except that the diagrams may contain components with
no vertices (closed circles). The relations $O_m$ and $P_{m+1}$ were
introduced by Le, Murakami and Ohtsuki
in~\cite{LeMurakamiOhtsuki:Universal}. $O_m$ says that disjoint union
with a closed circle is equivalent to multiplication by $(-2m)$ and
$P_{m+1}$ says that the sum of all ways of pairing up $2m+2$ stubs
attaching to the rest of the diagram is 0 (see Figure~\ref{fig:OP}).
By~\cite[Lemma~3.3]{LeMurakamiOhtsuki:Universal}, when we restrict to
the space $\Ao_{\leq m}(\emptyset)$ of diagrams of degree~$\leq m$, the
quotient by these relations is isomorphic to the corresponding
restriction $\calA_{\leq m}(\emptyset)$ of $\calA(\emptyset)$.  Hence
the assembly $\LMO_0$ can be regarded as taking values in
$\calA(\emptyset)$.

\begin{figure}[htpb]
\[ \printname{Orel}
        \setlength{\unitlength}{0.5\standardunitlength}
        \begin{array}{c}  \hspace{-1.7mm}
                \raisebox{-8pt}{\begingroup\makeatletter\ifx\SetFigFont\undefined%
\gdef\SetFigFont#1#2#3#4#5{%
  \reset@font\fontsize{#1}{#2pt}%
  \fontfamily{#3}\fontseries{#4}\fontshape{#5}%
  \selectfont}%
\fi\endgroup%
{\renewcommand{\dashlinestretch}{30}
\begin{picture}(4524,789)(0,-10)
\put(1212,387){\ellipse{600}{600}}
\path(1212,387)(1212,87)
\path(1212,387)(1472,537)
\path(1212,387)(952,537)
\put(4062,387){\ellipse{600}{600}}
\path(4062,387)(4062,87)
\path(4062,387)(4322,537)
\path(4062,387)(3802,537)
\path(3762,762)(3759,761)(3754,758)
	(3745,753)(3735,746)(3723,738)
	(3713,728)(3703,717)(3695,704)
	(3687,687)(3682,674)(3677,660)
	(3672,644)(3667,627)(3661,608)
	(3655,589)(3649,568)(3642,548)
	(3636,527)(3631,506)(3625,485)
	(3621,464)(3617,444)(3614,425)
	(3613,406)(3612,387)(3613,368)
	(3614,349)(3617,330)(3621,310)
	(3625,289)(3631,268)(3636,247)
	(3642,226)(3649,206)(3655,185)
	(3661,166)(3667,147)(3672,130)
	(3677,114)(3682,100)(3687,87)
	(3695,70)(3703,57)(3713,46)
	(3723,36)(3735,28)(3745,21)
	(3754,16)(3759,13)(3762,12)
\path(4362,762)(4365,761)(4370,758)
	(4379,753)(4390,746)(4401,738)
	(4411,728)(4421,717)(4429,704)
	(4437,687)(4442,674)(4447,660)
	(4452,644)(4457,627)(4463,608)
	(4469,589)(4475,568)(4482,548)
	(4488,527)(4493,506)(4499,485)
	(4503,464)(4507,444)(4510,425)
	(4511,406)(4512,387)(4511,368)
	(4510,349)(4507,330)(4503,310)
	(4499,289)(4493,268)(4488,247)
	(4482,226)(4475,206)(4469,185)
	(4463,166)(4457,147)(4452,130)
	(4447,114)(4442,100)(4437,87)
	(4429,70)(4421,57)(4411,46)
	(4401,36)(4390,28)(4379,21)
	(4370,16)(4365,13)(4362,12)
\put(462,387){\ellipse{600}{600}}
\path(162,762)(159,761)(154,758)
	(145,753)(135,746)(123,738)
	(113,728)(103,717)(95,704)
	(87,687)(82,674)(77,660)
	(72,644)(67,627)(61,608)
	(55,589)(49,568)(42,548)
	(36,527)(31,506)(25,485)
	(21,464)(17,444)(14,425)
	(13,406)(12,387)(13,368)
	(14,349)(17,330)(21,310)
	(25,289)(31,268)(36,247)
	(42,226)(49,206)(55,185)
	(61,166)(67,147)(72,130)
	(77,114)(82,100)(87,87)
	(95,70)(103,57)(113,46)
	(123,36)(135,28)(145,21)
	(154,16)(159,13)(162,12)
\path(1512,762)(1515,761)(1520,758)
	(1529,753)(1540,746)(1551,738)
	(1561,728)(1571,717)(1579,704)
	(1587,687)(1592,674)(1597,660)
	(1602,644)(1607,627)(1613,608)
	(1619,589)(1625,568)(1632,548)
	(1638,527)(1643,506)(1649,485)
	(1653,464)(1657,444)(1660,425)
	(1661,406)(1662,387)(1661,368)
	(1660,349)(1657,330)(1653,310)
	(1649,289)(1643,268)(1638,247)
	(1632,226)(1625,206)(1619,185)
	(1613,166)(1607,147)(1602,130)
	(1597,114)(1592,100)(1587,87)
	(1579,70)(1571,57)(1561,46)
	(1551,36)(1540,28)(1529,21)
	(1520,16)(1515,13)(1512,12)
\put(1812,312){\makebox(0,0)[lb]{\smash{{\mathmode{=(-2m)\cdot}}}}}
\end{picture}
} }
                \hspace{-1.9mm}
        \end{array}
 \qquad \printname{Prel}
        \setlength{\unitlength}{0.5\standardunitlength}
        \begin{array}{c}  \hspace{-1.7mm}
                \raisebox{-8pt}{\input draws/Prel.tex }
                \hspace{-1.9mm}
        \end{array}
 \]
\caption{
  The relation $O_m$ and the relation $P_2$. In the figure for $P_2$,
  the dashed square marks the parts of the diagrams where the relation
  is applied.
} \lbl{fig:OP}
\end{figure}

\end{myitemize}

A part of the divergence between the two definitions can be easily
remedied. In Section~\ref{sec:NDint-def} we will define a map
$\NDint:\calB(X)\to\Ao(\emptyset) / (O_m, P_{m+1})$ (called ``negative-dimensional formal integration'' for reasons to be explained in
Section~\ref{sec:integration}) and prove the following lemma:

\begin{lemma}\lbl{lem:equivalence}
  (Proof in Section~\ref{sec:NDint-def}, and see
  also~\cite[Lemma~6.3]{Le:Denominators})
  The composition $\NDint{}\circ\,\sigma$ is equal to $j_m/(O_m,P_{m+1})$.
\end{lemma}

Using this lemma we can redefine $\LMO_0^{(m)}$ to be the composition
$\NDint{}\circ\,\sigma\circ\Zcheck$. Comparing with
$\Aa_0=\Gint\circ\,\sigma\circ\Zcheck$ we see that the major difference
between $LMO_0$ and $\Aa_0$ is in the use of the different
``integrals'' $\NDint$ and $\Gint$. Thus the main technical challenge
in this paper is to compare the two integration theories. This is fully
achieved by the Proposition~\ref{prop:main} below, which says that
whenever $\Gint$ is defined, the two ``integrals'' differ only by a
normalization, and hence ultimately, the same holds for $\Aa$ and
$\LMO$.

\begin{definition} A ``function'' $G\in\calB(X)$ is said to be a perturbed
non-degenerate Gaussian (in the variables in $X$) if it is of the form
\[
  G=P\exp\Bigl(\frac12\sum_{x,y\in X}l_{xy}\strutn{x}{y}\Bigr)
\]
for some invertible symmetric matrix $\Lambda=(l_{xy})$ and some
$P\in\calB^+(X) \subset \calB(X)$. Here and throughout this paper, the
product on diagrams is the disjoint union product, $\strutn{x}{y}$
denotes a strut (a diagram in $\calB(X)$ made of a single edge with no
internal vertices) with ends labeled $x$ and $y$ and $\calB^+(X)$ is
the space of ``strutless'' diagrams in which each component has at
least one internal vertex
(cf.~\cite[Section~\ref{II-subsec:FGI}]{Aarhus:II}).
\end{definition}

\begin{proposition}\lbl{prop:main} (Proof in Section~\ref{sec:proof})
If $G$ is a perturbed non-degenerate Gaussian then
\[
  \NDint \!G\,dX = (-1)^{m|X|} (\det \Lambda)^m \Gint\!G\,dX
\]
in $\Ao_{\leq m}(\emptyset)/(O_m, P_{m+1})\simeq\calA_{\leq m}(\emptyset)$.
\end{proposition}

Note that $\NDint$ is defined in more cases than $\Gint$, but when they
are both defined, they are related in a simple way.

\subsection{Plan of paper}
In Section~\ref{sec:LMO}, we define $\NDint$, prove
Lemma~\ref{lem:equivalence}, and give an alternate formulation
of the $P_{m+1}$ relation.  In Section~\ref{sec:integration}
we prove some properties of $\NDint$ which justify the name
``negative-dimensional formal integration''.  These properties
are useful in Section~\ref{sec:proof}, where we prove the central
Proposition~\ref{prop:main} which is shown to imply Theorem~\ref{thm:main}
in the brief Section~\ref{sec:dull}.  Section~\ref{sec:philosophy}
has some philosophy on negative-dimensional spaces, sign choices for
diagrams, and the Rozansky-Witten invariants.
Appendix~\ref{sec:compexample} compares the definitions of $\Gint$ and
$\NDint$ by working out the first two terms in a non-trivial integral.


\section{A reformulation of the Le-Murakami-Ohtsuki invariant}\lbl{sec:LMO}

In this section we give a (minor) reformulation of the Le-Murakami-Ohtsuki
invariant $\LMO^{(m)}$. In Section~\ref{sec:NDint-def} we present our
definition of $\NDint$ and prove equivalence with the definition of
$j_m$ in~\cite{LeMurakamiOhtsuki:Universal}.  In Section~\ref{sec:Crel}
we state and prove an alternate form $C_{2m+1}$ of the relation $P_{m+1}$.

\subsection{Definition and notations}\lbl{sec:NDint-def}

\begin{definition} \lbl{def:NDint}
Let ``negative-dimensional integration''
 be defined by
\begin{equation} \lbl{eq:NDint-def}
  \begin{gathered}
  \NDint:\calB(X)\to\calA(\emptyset)/(O_m,P_{m+1})\\
  \NDint \!G\,dX = \Big\langle
  \prod_{x\in X} \frac{1}{m!} \left(\frac{\strutu{\prtl x}{\prtl x}}2
                                \right)^m \comma G
  \Big\rangle_{\!X} \Big/ (O_m, P_{m+1})
  \end{gathered}
\end{equation}
\end{definition}
Here the pairing $\langle\cdot\comma\cdot\rangle_X:
\calB(\partial_X)\otimes\calB(X)\to\Ao(\emptyset)$ is defined (like
in~\cite[Definition~\ref{I-def:Gint}]{Aarhus:I}) by
\[
\langle D_1\comma D_2\rangle_{X}=
  \left(\parbox{3.3in}{
    sum of all ways of gluing the $\prtl x$-marked legs of
    $D_1$ to the $x$-marked legs of $D_2$, for all $x\in X$
  }\right),
\]
where, as there, $\partial_X=\{\partial_x:x\in X\}$ denotes a set of
labels ``dual'' to the ones in $X$, and the sum is declared to be $0$
if the numbers of appropriately marked legs don't match.  If it is clear
which legs are to be attached, the subscript $X$ may be omitted.

In other words, $\NDint G$ is the composition of:
\begin{itemize}
\item projection of $G$ to the component with exactly $2m$ legs of
  each color in $X$
\item sum over all $\big((2m-1)!!\big)^{|X|} =
  \bigl(\frac{(2m)!}{2^m m!}\bigr)^{|X|}$
  ways of pairing up the legs of each color in $X$
\item quotient by the $O_m$ and $P_{m+1}$ relations.
\end{itemize}

An example of how to apply this definition is given in Appendix~\ref{sec:comp-ndint}.

Our definition of $\NDint$ is slightly different in appearance
than the definition of the corresponding object, $j_m$,
in~\cite{LeMurakamiOhtsuki:Universal}. For one, $j_m$ is defined on
$\calA(\uparrow_X)$ while $\NDint$ is defined on the different but
isomorphic space $\calB(X)$. We now prove Lemma~\ref{lem:equivalence},
which says that this is the only difference between the two maps.

\begin{proof}[Proof of Lemma~\ref{lem:equivalence}]

We prove that $j_m\circ\chi=\NDint$, where
$\chi:\calB(X)\to\calA(\uparrow_X)$ is the inverse of $\sigma$,
defined by mapping every $X$-marked uni-trivalent diagram in $\calB(X)$
to the average of all the ways of ordering and attaching its legs to
$n=|X|$ vertical arrows marked by the elements of $X$, respecting
the markings (\cite{Bar-Natan:Vassiliev, Bar-Natan:Homotopy},
\cite[Definition~\ref{I-def:sigma}]{Aarhus:I}).

First recall from~\cite{LeMurakamiOhtsuki:Universal} how $j_m$ is defined.
If $D$ is a diagram representing a class in $\calA(\uparrow_X)$, then
$j_m(D)$ is computed by removing the $n$ arrows from $D$ so that $n$
groups of stubs remain, and then by gluing certain (linear combinations
of) forests on these stubs, so that each tree in each forest gets glued
only to the stubs within some specific group. It is not obvious that $j_m$
is well defined; it may not respect the $STU$ relation. With some effort,
it is proven in~\cite{LeMurakamiOhtsuki:Universal} that for the specific
combinations of forests used there, $j_m$ is indeed well defined.

Now every tree that has internal vertices has some two leafs that connect
to the same internal vertex, and hence (modulo $AS$), every such tree is
anti-symmetric modulo some transposition of its leafs. Thus gluing such a
tree to symmetric combinations of diagrams, such as the ones in the image
of $\chi$, we always get $0$. Hence in the computation of $j_m\circ\chi$
it is enough to consider forests of trees that have no internal vertices;
that is, forests of struts. Extracting the precise coefficients
from~\cite{LeMurakamiOhtsuki:Universal} one easily sees that they are
the same as in~\eqref{eq:NDint-def}, and hence $j_m\circ\chi=\NDint$,
as required.  
\end{proof}

\subsection{The $C_l$ relations} \lbl{sec:Crel}
It will be convenient, for use in Proposition~\ref{prop:transinv}, to
give another reformulation of the definitions
of~\cite{LeMurakamiOhtsuki:Universal}.  Instead of their $P_{m+1}$
relation, we may use another relation, the $C_{2m+1}$ relation.  For
motivation for this relation, see Section~\ref{sec:OCrels}.

\begin{definition} The $C_l$ relation\footnote{`$C$' for `Crocodile'.}
applies when we have a diagram with two sets of $l$ stubs (or teeth)
each, and says that the sum of the diagrams obtained by attaching the
two sets of stubs to each other in all $l!$ possible ways is 0, as in
the following diagram.
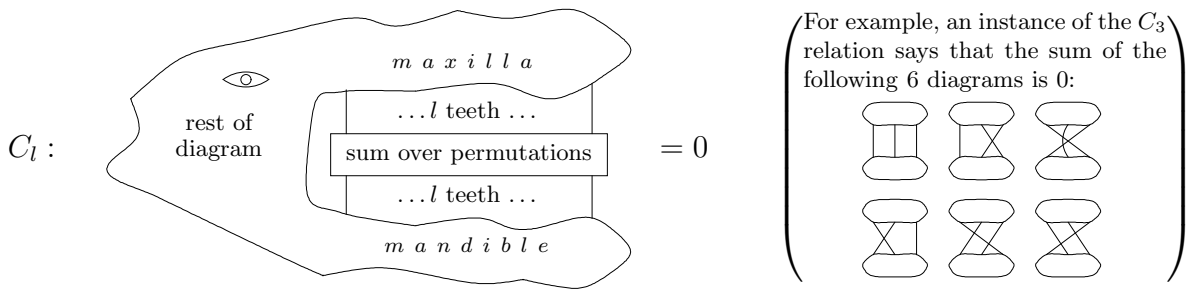
\[
  {
    \def\maxilla{{\scriptstyle m\ a\ x\ i\ l\ l\ a}}
    \def\mandible{{\scriptstyle m\ a\ n\ d\ i\ b\ l\ e}}
    C_l:\quad\printname{Crel}
        \setlength{\unitlength}{0.32\standardunitlength}
        \begin{array}{c}  \hspace{-1.7mm}
                \raisebox{-8pt}{\input draws/Crel.tex }
                \hspace{-1.9mm}
        \end{array}
 = 0
  }\qquad\left(\parbox{1.9in}{\scriptsize
    For example, an instance of the $C_3$ relation says that the sum of
    the following $6$ diagrams is $0$:
    \begin{center} $\printname{C3}
        \setlength{\unitlength}{0.45\standardunitlength}
        \begin{array}{c}  \hspace{-1.7mm}
                \raisebox{-8pt}{\input draws/C3.tex }
                \hspace{-1.9mm}
        \end{array}
$ \end{center}
  }\right)
\]
\end{definition}

Note that both sets of relations, the $P_m$'s and the $C_l$'s, are
decreasing in power. Namely, $P_m$ implies $P_{m+1}$ and $C_l$ implies
$C_{l+1}$, for every $l$ and $m$ (one easily sees that $P_{m+1}$ is a sum
of instances of $P_m$, and likewise for $C_{l+1}$ and $C_l$). The lemma
below says that up to an index-doubling, the two chains of relations are
equivalent.

\begin{lemma}\lbl{lem:Crel}
The relations $C_{2m+1}$, $C_{2m+2}$, and $P_{m+1}$ are equivalent.
(All of these relations may be applied inside any space of diagrams,
regardless of the IHX, STU, or any other relations, so long as closed
circles are allowed).
\end{lemma}

\par\noindent
\parbox{5.2in}{
  {\em Proof.}
  It was already noted that $C_{2m+1}$ implies $C_{2m+2}$. Next, it is
  easy to see that $C_{2m+2}$ implies $P_{m+1}$: just apply the relation
  $C_{2m+2}$ in the diagram shown on the right, and you get (a positive
  multiple of) the relation $P_{m+1}$.
}\ \hfill$\printname{CtoP}
        \setlength{\unitlength}{0.4\standardunitlength}
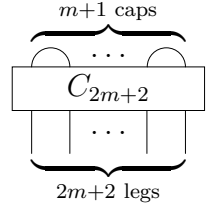
        \begin{array}{c}  \hspace{-1.7mm}
                \raisebox{-8pt}{%
\begingroup\makeatletter\ifx\SetFigFont\undefined
\def\x#1#2#3#4#5#6#7\relax{\def\x{#1#2#3#4#5#6}}%
\expandafter\x\fmtname xxxxxx\relax \def\y{splain}%
\ifx\x\y   
\gdef\SetFigFont#1#2#3{%
  \ifnum #1<17\tiny\else \ifnum #1<20\small\else
  \ifnum #1<24\normalsize\else \ifnum #1<29\large\else
  \ifnum #1<34\Large\else \ifnum #1<41\LARGE\else
     \huge\fi\fi\fi\fi\fi\fi
  \csname #3\endcsname}%
\else
\gdef\SetFigFont#1#2#3{\begingroup
  \count@#1\relax \ifnum 25<\count@\count@25\fi
  \def\x{\endgroup\@setsize\SetFigFont{#2pt}}%
  \expandafter\x
    \csname \romannumeral\the\count@ pt\expandafter\endcsname
    \csname @\romannumeral\the\count@ pt\endcsname
  \csname #3\endcsname}%
\fi
\fi\endgroup
{\renewcommand{\dashlinestretch}{30}
\begin{picture}(3024,3675)(0,-10)
\put(2412.000,2250.000){\arc{600.000}{3.1416}{6.2832}}
\put(612.000,2250.000){\arc{600.000}{3.1416}{6.2832}}
\path(12,2250)(3012,2250)(3012,1575)
	(12,1575)(12,2250)
\path(2113,1574)(2113,899)
\path(2713,1574)(2713,899)
\path(913,1574)(913,899)
\path(313,1574)(313,899)
\put(1512,2400){\makebox(0,0)[b]{\smash{{\mathmode{\ldots}}}}}
\put(1513,1199){\makebox(0,0)[b]{\smash{{\mathmode{\ldots}}}}}
\put(1437,3525){\makebox(0,0)[lb]{\smash{{\mathmode{\relax}}}}}
\put(312,2625){\makebox(0,0)[lb]{\smash{{\mathmode{\overbrace{\hskip.8in}^{m+1 \text{ caps}}}}}}}
\put(313,824){\makebox(0,0)[lb]{\smash{{\mathmode{\underbrace{\hskip.8in}_{2m+2 \text{ legs}}}}}}}
\put(1512,1800){\makebox(0,0)[b]{\smash{{\mathmode{\textstyle{C_{2m+2}}}}}}}
\put(1212,0){\makebox(0,0)[lb]{\smash{{\mathmode{\relax}}}}}
\end{picture}
}
 }
                \hspace{-1.9mm}
        \end{array}
$\hfill\ 

\par\ \par\noindent
\ \hfill${
  \def\gluings{{\scriptstyle g\ l\ u\ i\ n\ g\ s}}
  \printname{caps}
        \setlength{\unitlength}{0.4\standardunitlength}
        \begin{array}{c}  \hspace{-1.7mm}
                \raisebox{-8pt}{\begingroup\makeatletter\ifx\SetFigFont\undefined%
\gdef\SetFigFont#1#2#3#4#5{%
  \reset@font\fontsize{#1}{#2pt}%
  \fontfamily{#3}\fontseries{#4}\fontshape{#5}%
  \selectfont}%
\fi\endgroup%
{\renewcommand{\dashlinestretch}{30}
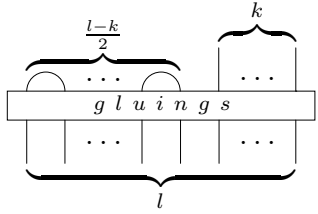
\begin{picture}(4824,3226)(0,-10)
\put(2412.000,1875.000){\arc{600.000}{3.1416}{6.2832}}
\put(612.000,1875.000){\arc{600.000}{3.1416}{6.2832}}
\path(3312,1875)(3312,2551)
\drawline(4512,1875)(4512,1875)
\path(4512,1875)(4512,2551)
\path(12,1875)(4812,1875)(4812,1425)
	(12,1425)(12,1875)
\path(2113,1424)(2113,749)
\path(2713,1424)(2713,749)
\path(3313,1424)(3313,749)
\path(4513,1424)(4513,749)
\path(913,1424)(913,749)
\path(313,1424)(313,749)
\put(1512,2025){\makebox(0,0)[b]{\smash{{\mathmode{\ldots}}}}}
\put(1513,1049){\makebox(0,0)[b]{\smash{{\mathmode{\ldots}}}}}
\put(3913,1049){\makebox(0,0)[b]{\smash{{\mathmode{\ldots}}}}}
\put(2412,0){\makebox(0,0)[lb]{\smash{{\mathmode{\relax}}}}}
\put(3912,2040){\makebox(0,0)[b]{\smash{{\mathmode{\ldots}}}}}
\put(313,674){\makebox(0,0)[lb]{\smash{{\mathmode{\underbrace{\hskip1.4in}_{l}}}}}}
\put(311,2176){\makebox(0,0)[lb]{\smash{{\mathmode{\overbrace{\hskip.8in}^{\frac{l-k}{2}}}}}}}
\put(611,3076){\makebox(0,0)[lb]{\smash{{\mathmode{\relax}}}}}
\put(3311,2626){\makebox(0,0)[lb]{\smash{{\mathmode{\overbrace{\hskip.4in}^k}}}}}
\put(1361,1576){\makebox(0,0)[lb]{\smash{{\mathmode{\gluings}}}}}
\end{picture}
} }
                \hspace{-1.9mm}
        \end{array}
}
$\hfill\ 
\parbox{4.6in}{
  The proof of the implication $P_{m+1}\Rightarrow
  C_{2m+1}$ is essentially the proof of Lemma~3.1
  of~\cite{LeMurakamiOhtsuki:Universal}, though the result is different.
  First, a definition: for $k\leq l$ and $l-k$ even, the diagram part
  $C^k_l$ has $k$ legs pointing up and $l$ legs pointing down, and it
  is the sum of all ways of attaching all of the $k$ legs to some of the
  $l$ legs and then pairing up the remaining $(l-k)$ legs, as illustrated
  on the left.
}

\par\ \par
We now prove by induction that for all $0\leq k\leq m$, $C^{2k+1}_{2m+1}=0$
modulo $P_{m+1}$. For $k = 0$, $C^1_{2m+1}$ is just a version of
$P_{m+1}$.  For $k > 0$, apply $P_{m+k+1}$ (a consequence of $P_{m+1}$)
to a diagram with $2m+1$ legs pointing down and $2k+1$ legs pointing
up, like this:
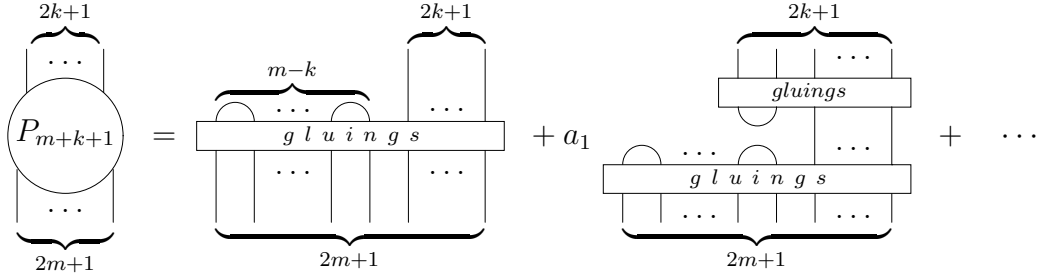
\[
  \def\glu{{\scriptstyle gluings}}
  \def\gluings{{\scriptstyle g\ l\ u\ i\ n\ g\ s}}
  \printname{Pmk1}
        \setlength{\unitlength}{0.4\standardunitlength}
        \begin{array}{c}  \hspace{-1.7mm}
                \raisebox{-8pt}{%
\begingroup\makeatletter\ifx\SetFigFont\undefined
\def\x#1#2#3#4#5#6#7\relax{\def\x{#1#2#3#4#5#6}}%
\expandafter\x\fmtname xxxxxx\relax \def\y{splain}%
\ifx\x\y   
\gdef\SetFigFont#1#2#3{%
  \ifnum #1<17\tiny\else \ifnum #1<20\small\else
  \ifnum #1<24\normalsize\else \ifnum #1<29\large\else
  \ifnum #1<34\Large\else \ifnum #1<41\LARGE\else
     \huge\fi\fi\fi\fi\fi\fi
  \csname #3\endcsname}%
\else
\gdef\SetFigFont#1#2#3{\begingroup
  \count@#1\relax \ifnum 25<\count@\count@25\fi
  \def\x{\endgroup\@setsize\SetFigFont{#2pt}}%
  \expandafter\x
    \csname \romannumeral\the\count@ pt\expandafter\endcsname
    \csname @\romannumeral\the\count@ pt\endcsname
  \csname #3\endcsname}%
\fi
\fi\endgroup
{\renewcommand{\dashlinestretch}{30}
\begin{picture}(1822,3750)(0,-10)
\path(1511,3225)(1511,2550)
\path(311,3225)(311,2550)
\path(1661,1350)(1661,525)
\path(161,1350)(161,525)
\put(311,3600){\makebox(0,0)[lb]{\smash{{\mathmode{\relax}}}}}
\put(161,0){\makebox(0,0)[lb]{\smash{{\mathmode{\relax}}}}}
\put(911,1875){\ellipse{1806}{1806}}
\put(911,1800){\makebox(0,0)[b]{\smash{{\mathmode{\textstyle P_{m+k+1}}}}}}
\put(911,675){\makebox(0,0)[b]{\smash{{\mathmode{\ldots}}}}}
\put(911,3000){\makebox(0,0)[b]{\smash{{\mathmode{\ldots}}}}}
\put(161,375){\makebox(0,0)[lb]{\smash{{\mathmode{\underbrace{\hskip.5in}_{2m+1}}}}}}
\put(386,3300){\makebox(0,0)[lb]{\smash{{\mathmode{\overbrace{\hskip.35in}^{2k+1}}}}}}
\end{picture}
}
 }
                \hspace{-1.9mm}
        \end{array}
 = \printname{caps2}
        \setlength{\unitlength}{0.4\standardunitlength}
        \begin{array}{c}  \hspace{-1.7mm}
                \raisebox{-8pt}{%
\begingroup\makeatletter\ifx\SetFigFont\undefined
\def\x#1#2#3#4#5#6#7\relax{\def\x{#1#2#3#4#5#6}}%
\expandafter\x\fmtname xxxxxx\relax \def\y{splain}%
\ifx\x\y   
\gdef\SetFigFont#1#2#3{%
  \ifnum #1<17\tiny\else \ifnum #1<20\small\else
  \ifnum #1<24\normalsize\else \ifnum #1<29\large\else
  \ifnum #1<34\Large\else \ifnum #1<41\LARGE\else
     \huge\fi\fi\fi\fi\fi\fi
  \csname #3\endcsname}%
\else
\gdef\SetFigFont#1#2#3{\begingroup
  \count@#1\relax \ifnum 25<\count@\count@25\fi
  \def\x{\endgroup\@setsize\SetFigFont{#2pt}}%
  \expandafter\x
    \csname \romannumeral\the\count@ pt\expandafter\endcsname
    \csname @\romannumeral\the\count@ pt\endcsname
  \csname #3\endcsname}%
\fi
\fi\endgroup
{\renewcommand{\dashlinestretch}{30}
\begin{picture}(4824,3750)(0,-10)
\put(3311,3600){\makebox(0,0)[lb]{\smash{{\mathmode{\relax}}}}}
\put(311,0){\makebox(0,0)[lb]{\smash{{\mathmode{\relax}}}}}
\put(612.000,2099.000){\arc{600.000}{3.1416}{6.2832}}
\path(3312,2099)(3312,3225)
\drawline(4512,2099)(4512,2099)
\path(4512,2099)(4512,3225)
\path(12,2099)(4812,2099)(4812,1649)
	(12,1649)(12,2099)
\path(2113,1648)(2113,525)
\path(2713,1648)(2711,525)
\path(3313,1648)(3313,525)
\path(4513,1648)(4513,525)
\path(913,1648)(913,525)
\put(1512,2249){\makebox(0,0)[b]{\smash{{\mathmode{\ldots}}}}}
\put(1513,1273){\makebox(0,0)[b]{\smash{{\mathmode{\ldots}}}}}
\put(3913,1273){\makebox(0,0)[b]{\smash{{\mathmode{\ldots}}}}}
\put(3912,2264){\makebox(0,0)[b]{\smash{{\mathmode{\ldots}}}}}
\put(311,2400){\makebox(0,0)[lb]{\smash{{\mathmode{\overbrace{\hskip.8in}^{m-k}}}}}}
\put(1361,1800){\makebox(0,0)[lb]{\smash{{\mathmode{\gluings}}}}}
\put(3311,3300){\makebox(0,0)[lb]{\smash{{\mathmode{\overbrace{\hskip.4in}^{2k+1}}}}}}
\put(311,450){\makebox(0,0)[lb]{\smash{{\mathmode{\underbrace{\hskip1.4in}_{2m+1}}}}}}
\path(313,1648)(313,525)
\put(2412.000,2099.000){\arc{600.000}{3.1416}{6.2832}}
\end{picture}
}
 }
                \hspace{-1.9mm}
        \end{array}
 + a_1
  \printname{MoreCaps}
        \setlength{\unitlength}{0.4\standardunitlength}
        \begin{array}{c}  \hspace{-1.7mm}
                \raisebox{-8pt}{\begingroup\makeatletter\ifx\SetFigFont\undefined
\def\x#1#2#3#4#5#6#7\relax{\def\x{#1#2#3#4#5#6}}%
\expandafter\x\fmtname xxxxxx\relax \def\y{splain}%
\ifx\x\y   
\gdef\SetFigFont#1#2#3{%
  \ifnum #1<17\tiny\else \ifnum #1<20\small\else
  \ifnum #1<24\normalsize\else \ifnum #1<29\large\else
  \ifnum #1<34\Large\else \ifnum #1<41\LARGE\else
     \huge\fi\fi\fi\fi\fi\fi
  \csname #3\endcsname}%
\else
\gdef\SetFigFont#1#2#3{\begingroup
  \count@#1\relax \ifnum 25<\count@\count@25\fi
  \def\x{\endgroup\@setsize\SetFigFont{#2pt}}%
  \expandafter\x
    \csname \romannumeral\the\count@ pt\expandafter\endcsname
    \csname @\romannumeral\the\count@ pt\endcsname
  \csname #3\endcsname}%
\fi
\fi\endgroup
{\renewcommand{\dashlinestretch}{30}
\begin{picture}(4824,3750)(0,-10)
\put(312,0){\makebox(0,0)[lb]{\smash{{\mathmode{\relax}}}}}
\put(2112,3600){\makebox(0,0)[lb]{\smash{{\mathmode{\relax}}}}}
\put(2412.000,2323.000){\arc{600.000}{6.2832}{9.4248}}
\path(3312,1424)(3312,2325)
\drawline(4512,1424)(4512,1424)
\path(4512,1424)(4512,2325)
\path(12,1424)(4812,1424)(4812,974)
	(12,974)(12,1424)
\path(2113,973)(2113,525)
\path(2713,973)(2713,525)
\path(3313,973)(3313,525)
\path(4513,973)(4513,525)
\path(913,973)(913,525)
\path(313,973)(313,525)
\path(4812,2323)(1812,2323)(1812,2773)
	(4812,2773)(4812,2323)
\path(4512,2773)(4512,3225)
\path(3312,2773)(3312,3225)
\path(2712,2773)(2712,3225)
\path(2112,2773)(2112,3225)
\put(1512,1574){\makebox(0,0)[b]{\smash{{\mathmode{\ldots}}}}}
\put(1513,598){\makebox(0,0)[b]{\smash{{\mathmode{\ldots}}}}}
\put(3913,598){\makebox(0,0)[b]{\smash{{\mathmode{\ldots}}}}}
\put(3912,1648){\makebox(0,0)[b]{\smash{{\mathmode{\ldots}}}}}
\put(312,450){\makebox(0,0)[lb]{\smash{{\mathmode{\underbrace{\hskip1.4in}_{2m+1}}}}}}
\put(1362,1125){\makebox(0,0)[lb]{\smash{{\mathmode{\gluings}}}}}
\put(3912,3073){\makebox(0,0)[b]{\smash{{\mathmode{\ldots}}}}}
\put(2112,3300){\makebox(0,0)[lb]{\smash{{\mathmode{\overbrace{\hskip.8in}^{2k+1}}}}}}
\put(2637,2475){\makebox(0,0)[lb]{\smash{{\mathmode{\glu}}}}}
\put(612.000,1424.000){\arc{600.000}{3.1416}{6.2832}}
\put(2412.000,1424.000){\arc{600.000}{3.1416}{6.2832}}
\end{picture}
} }
                \hspace{-1.9mm}
        \end{array}
 + \quad\cdots
\]
As shown in the diagram, the result splits into a sum over the number
$2l$ of the upwards pointing legs that get paired with each other;
for $l=0$, we just get $C^{2k+1}_{2m+1}$; for $l > 0$, the result
can be considered to split into two diagrams, a (positive multiple of
a) reversed $C^{2(k-l)+1}_{2k+1}$ on top of a $C^{2(k-l)+1}_{2m+1}$.
But, by the induction hypothesis, this latter term $C^{2(k-l)+1}_{2m+1}$
vanishes modulo $P_{m+1}$, and we are left with just the first term,
which is therefore also a consequence of $P_{m+1}$. This completes the
inductive proof.  To conclude the proof of Lemma~\ref{lem:Crel}, note
that $C_{2m+1} = C{}^{2m+1}_{2m+1}$. \qed


\section{Negative-dimensional formal integration} \lbl{sec:integration}

In this section, we give several justifications of the name
``negative-dimensional formal integration'' for the map $\NDint$
defined above.  While doing this, we prove several properties of $\NDint$
(Propositions~\ref{prop:transinv} and~\ref{prop:gaussint}) that are used
in the proof of Proposition~\ref{prop:main} in Section~\ref{sec:proof}.

\subsection{Why integration?}

First, why should $\NDint$ be called an integral?  In general, an integral
is (more or less) a linear map from some space of functions to the corresponding
space of scalars. In our case, the appropriate space of ``functions'' is
$\calB(X)$ and the appropriate space of ``scalars'' is $\calA(\emptyset)$.%
\footnote{%
More on the interpretation of diagrams as functions and/or scalars
appears in~\cite[Section~\ref{I-subsec:DiagrammaticCase}]{Aarhus:I}
and~\cite[Section~\ref{II-sec:calculus}]{Aarhus:II}.} The linearity of
$\NDint$ is immediate.

But $\NDint$ is not just any integral; it is a Lebesgue integral.
The defining property of the usual Lebesgue integral on $\RR^n$ is
translation invariance by a vector $(\bar x^i)$. We show that the parallel
property holds for $\NDint$:

\begin{proposition}[Translation Invariance]\lbl{prop:transinv}
For any diagram $D \in \calB(X)$, we have 
\begin{equation}\lbl{eq:transinv}
  \NDint\!D\,dX = \NDint\!D / (x\mapsto x + \bar x) \, dX
\end{equation}
\end{proposition}

The notation $D /(x\mapsto x + \bar x)$ means (as
in~\cite[Section~\ref{II-subsec:GeneralSetup}]{Aarhus:II}), for each
leg of $D$ colored $x$ for $x\in X$, sum over coloring the leg by $x$
or by $\bar x$.  (So we end up with a sum of $2^t$ terms, where $t$ is
the number of $X$ colored legs in $D$.)  The set $\bar{X}=\{\bar x \mid
x \in X\}$ is an independent set of variables for the formal translation.

\begin{proof}
  By the relation $P_{m+1}$ in the definition of $\NDint$ (or rather,
by $C_{2m+1}$), any diagram $D\in\calB(X)$ with more than $2m$ legs on
any component gets mapped to 0 on either side of~(\ref{eq:transinv}),
so we may assume that $D$ has $2m$ legs or fewer of any color.  But,
for the right hand side to be non-zero, $D/(x\mapsto x+\bar x)$ must
have exactly $2m$ legs colored $x$ for each $x \in X$; this can only
happen from diagrams in $D$ with exactly $2m$ legs of each color and
when none of them get converted to $\bar x$.  But these are exactly
the diagrams appearing in the integral on the left hand side.
\end{proof}

\subsection{The relations $O_m$ and $C_{2m+1}$}\lbl{sec:OCrels}
Now that you're convinced that $\NDint$ is an integral, you must be
wondering why we called it a ``negative-dimensional'' integral.
Recall what $\NDint$ is: it is the sum over all ways of gluing in some
struts, followed by the quotient by the $O_m$ and $P_{m+1}$
relations.  This quotient is crucial; otherwise $\NDint$ is some
random map without particularly nice properties.  But what are these
relations?

The relation $O_m$ is simple.
Recall~\cite[Section~\ref{II-sec:calculus}]{Aarhus:II} that we like to
think of diagrams as representing tensors and/or functions in/on some
vector space $V$. Since a strut corresponds to the identity tensor in
$V^\star\otimes V$ (cf.~\cite[Figure~\ref{II-fig:StrutLaws}]{Aarhus:II}),
its closure, a circle, should correspond to the trace of the identity, or
the dimension of $V$. Hence $O_m$, which says that a circle is equivalent
to the constant $(-2m)$, is the parallel of saying ``$\dim V=-2m$''.

The relation $P_{m+1}$ is a little more subtle. It is easier to look
at the equivalent relation $C_{2m+1}$, which implies the relation $C_l$
for every $l>2m$. If a single vertex corresponds to some space $V$, then
a collection of $l$ vertices corresponds to $V^{\otimes l}$; and, when
we sum over all permutations without signs, we get (a multiple of) the
projection onto the symmetric subspace, $S^l(V)$. The relation $C_l$ says
that this projection (and hence the target, $S^l(V)$) is $0$. Compare
this with the following statement about $\RR^k$ for $k\geq 0$:
\[
  \dim S^l(\RR^k)
  = \left(\begin{array}{c} l+k-1 \\ l \end{array}\right)
  = \frac{(k+l-1)(k+l-2)\cdots(k+1)k}{l(l-1)\cdots 2 \cdot 1}.
\]
We can see from this formula that if a space $V$ formally has a dimension
$k=-2m$, then $\dim S^l(V)$ vanishes precisely when $l>2m$. This is in
complete agreement with what we just found about $P_{m+1}$.

\subsection{An example: Gaussian integration} \lbl{subsec:compute}

Enough of generalities, let's compute! Consider the well known Gaussian
integral over $\RR^n$,
\[
  \int_{\RR^n} e^{q(x,x)/2}\,d^n x =
  \frac{(2\pi)^{\frac{n}{2}}}{(\det\,-q)^{\frac12}}
\]
where $q$ is an arbitrary negative-definite quadratic form $q$.
The factor $(2\pi)^{n/2}$ is just a normalization factor that could
be absorbed into the measure $d^n x$.  (Recall that we identified
$\NDint$ as Lebesgue integration by translation invariance, which only
determines the measure up to an overall scale factor.)  The remaining
factor, $(\det\,-q)^{-1/2}$, is the more fundamental one.  Does a similar
result hold for $\NDint$?  In order to answer this question, one would
first have to know what a ``determinant'' of a quadratic form on a
negative-dimensional space is.  While there is a good answer to this
question (called the ``superdeterminant'' or the ``Berezinian''; see,
e.g., \cite[page~82]{Berezin:Superanalysis}), it would take us too
far afield to discuss it in full.  Instead, let us take a slightly
different tack.  Fix a negative-definite quadratic form $\Lambda$
on $\RR^n$ once and for all, and consider the quadratic form $q =
\Lambda\otimes(\delta_{ij})$ on $\RR^n\otimes\RR^k\cong\RR^{nk}$. We
have $\det q=(\det\Lambda)^k$ and so we find that
\begin{equation}\lbl{eq:pdGaussInt}
  \int_{\RR^{nk}} e^{q(x,x)/2}\,d^{nk}x = C (\det\,-\Lambda)^{-k/2}
\end{equation}
for some constant $C$.

Consider now the sum $\sum_{x,y\in X}l_{xy}\strutn{x}{y}$
in $\calB(X)$. According to the voodoo of diagrammatic
calculus~\cite[Section~\ref{II-sec:calculus}]{Aarhus:II}, it is in
analogy with a quadratic form on $V^{\otimes n}$, where $n=|X|$ and
$V$ is some vector space that plays a role similar to $\RR^k$ in the
above discussion. The proposition below is then the diagrammatic analog
of~\eqref{eq:pdGaussInt}, taking $k=\dim V=-2m$.

\begin{proposition}\lbl{prop:gaussint}
For any set $X$ with $|X|=n$ and $\Lambda=(l_{xy})$ a symmetric matrix
on $\RR^X$,
\[ \NDint\!\exp\left(\frac12 \sum_{x,y\in X} l_{xy} \strutn{x}{y}\right)
  = (\det\, -\Lambda)^m = (-1)^{nm}(\det \Lambda)^m
\]
\end{proposition}

Note that other than symmetry, there is no restriction on the
matrix $\Lambda$. This proposition is a consequence of Lemma~4.2
of~\cite{LeMurakamiOhtsuki:Universal} and a computation for $m=1$
(given in~\cite{LeMurakami2Ohtsuki:3Manifold}), but we give our own
proof for completeness, and also to provide a more direct link to typical
determinant calculations.

\begin{proof}
We are to calculate the reduction modulo $O_m$ and $P_{m+1}$ of
\[ D_1 := \Bigg\langle
    \prod_{x\in X} \frac{1}{m!}
     \left(\frac{\strutu{\prtl x}{\prtl x}}{2}\right)^{\!m}\!\comma 
       \exp\left(\frac12 \sum_{x,y\in X} l_{xy}\strutn{x}{y}\right)
  \Bigg\rangle_X.
\]
The only terms that can appear in $D_1$ are closed loops. The relation
$O_m$ replaces each of these by a number, reducing the result to
$\QQ$. The relation $P_{m+1}$ is irrelevant and will be ignored in the
remainder of the proof.

Introduce a new set of variables $A$ (and dual variables $\partial_A$)
with $|A| = m$, and consider
\begin{equation} \lbl{eq:D2}
   D_2 := \Bigg\langle
    \prod_{x\in X} \frac{1}{m!} \left(\sum_{a\in A}
      \downstrut{(x,a)}{(\prtl x,\prtl a)}\right)^{\!m}\!\comma
    \exp\left(\sum_{x,y\in X}\sum_{a\in A} l_{xy}
      \upstrut{(y, a)}{(\prtl x, \prtl a)}\right)
  \Bigg\rangle_{X A}.
\end{equation}

Let us compare $D_1$ and $D_2$. After all relevant gluings, $D_1$
becomes a sum (with coefficients) of disjoint unions of unoriented
loops, each of which is a polygon of struts whose vertices are colored
by elements of $X$. Similarly, $D_2$ is also a sum (with coefficients)
of a disjoint union of loops, only that now the loops are oriented and
the struts they are made of are colored by the elements of $X A$,
keeping the $A$ part of the coloring constant along each loop. In both
cases the coefficients come from the same simple rule, which involves
only the $X$ part of the coloring. We see that each term in $D_1$ with
$c$ circles corresponds to $(2m)^c$ terms of $D_2$: for each loop in a
given term of $D_1$, choose a color $a\in A$ and an orientation, and
you get a term in $D_2$. So we find that
\[ D_2 \big/ (\bigcirc = -1)\,=\,D_1 \big/ (\bigcirc = -2m). \]

Recall that $\frac{1}{4!}(a+b+c+d)^4=abcd+(\text{non-multi-linear terms})$.
Similarly,
\[
  \prod_{x\in X} \frac{1}{m!} \left(
    \sum_{a\in A} \downstrut{(x,a)}{(\prtl x,\prtl a)}
  \right)^{\!m}
  = \prod_{(x,a)\in X  A} \downstrut{(x, a)}{(\prtl x, \prtl a)}
  + (\text{terms with strut repetitions}).
\]
We assert that terms with strut repetitions can be ignored
in the computation of $D_2\big/(\bigcirc=-1)$. Indeed, for
some fixed $x_0\in X$ and $a_0\in A$ set $\alpha=(x_0,a_0)$ and
$\prtl\alpha=(\prtl{x_0},\prtl{a_0})$, and suppose a strut repetition
like $\downstrut{\alpha}{\prtl\alpha} \downstrut{\alpha}{\prtl\alpha}$
occurs within the left operand of a pairing as in~\eqref{eq:D2}. Then,
as illustrated in Figure~\ref{fig:StrutRepetition}, the gluings in
the evaluation of the pairings come in pairs. One easily sees that the
number of cycles differs by $ 1$ for the gluings within each pair,
and hence modulo $(\bigcirc=-1)$ the whole sum of gluings vanishes.

\begin{figure}[htpb]
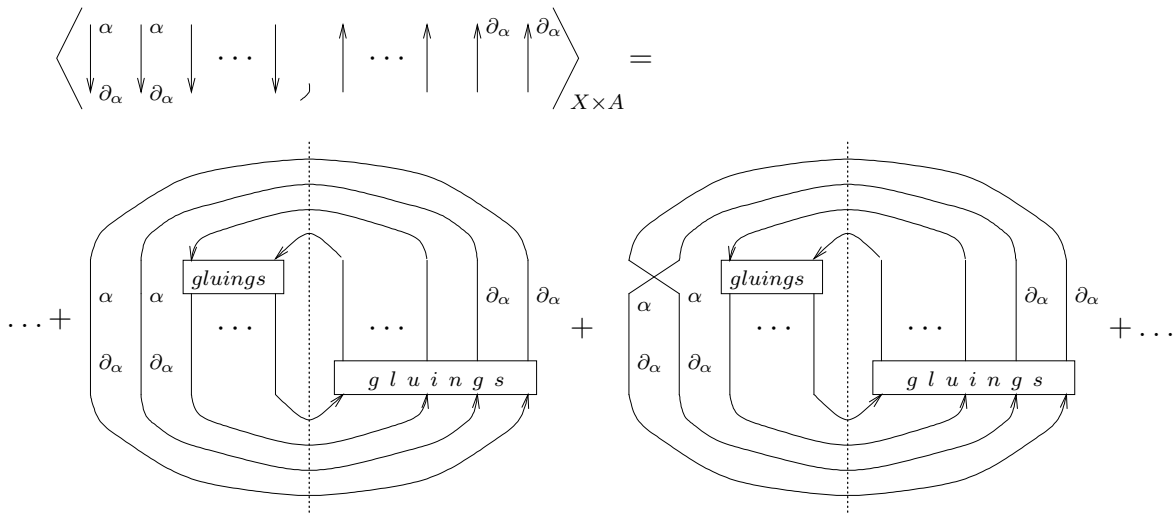

\[
  \def\a{\scriptstyle\alpha}
  \def\pa{\scriptstyle\prtl\alpha}
  \def\glu{{\scriptstyle gluings}}
  \def\gluings{{\scriptstyle g\ l\ u\ i\ n\ g\ s}}
  \printname{StrutRepetition}
        \setlength{\unitlength}{0.7\standardunitlength}
        \begin{array}{c}  \hspace{-1.7mm}
                \raisebox{-8pt}{\input draws/StrutRepetition.tex }
                \hspace{-1.9mm}
        \end{array}

\]
\caption{Two ways of gluing a repeating strut.} \lbl{fig:StrutRepetition}
\end{figure}

Now compare~$D_2$ to
\[ D_3 :=
  \Bigg\langle \prod_{(x,a)\in X  A}
    \downstrut{(x, a)}{(\prtl x, \prtl a)} \comma
    \exp\left(\sum_{x,y\in X}\sum_{a\in A} l_{xy}
      \upstrut{(y,a)}{(\prtl x,\prtl a)}\right)
  \Bigg\rangle_{X A}.
\]
By our assertion, $D_2=D_3$ modulo $(\bigcirc=-1)$. But $D_3$ looks very
much like the usual formula for the determinant: it reduces to
\[
  D_3 \big/ (\bigcirc=-1)\,= \sum_{\pi\in S(X A)}
    \prod_{(x,a)\in X A} l_{xa,\,\pi(xa)}(-1)^{\text{cycles}(\pi)}
\]
where $(l_{xa,yb})=\Lambda\otimes(\delta_{ab})$. Using the relationship
between the number of cycles of a permutation $\pi\in S(X A)$
and its signature, $(-1)^{\text{cycles}(\pi)}=(-1)^{nm}\sgn(\pi)$,
we find that
\[ \NDint G = (-1)^{nm}\det (l_{xa,yb}) = (-1)^{nm}(\det \Lambda)^m \]
as required.
\end{proof}


\section{Relating the two integration theories}\lbl{sec:proof}

The classical computation of perturbed Gaussian integration
uses only translation invariance, and a single non-perturbed
computation to determine the normalization coefficient. For negative
dimensional integration, translation invariance was proven in
Proposition~\ref{prop:transinv}, and the non-perturbed computation is in
Section~\ref{subsec:compute}. So the proof of Proposition~\ref{prop:main}
proceeds just as in the classical computation of perturbed Gaussian
integration:

\begin{proof}[Proof of Proposition~\ref{prop:main}]
Recall that we are comparing $\NDint$ to $\Gint$ on a 
perturbed Gaussian with variables in $X$ and a non-degenerate
quadratic part $\Lambda = (l_{xy})$:
\[
  G=P\exp\biggl(\frac12 \,l_{xy}\strutn{x}{y}\biggr).
\]
(Here and throughout this proof,
repeated variables should be summed over $X$.)
From~\cite[Definition~\ref{I-def:Gint}]{Aarhus:I} we have
\[ \Gint\!G\,dX
  = \left\langle
    \exp\biggl(-\frac12\,l^{xy}\strutu{\prtl x}{\prtl y}\biggr), P
  \right\rangle,
\]
with, as usual, $(l^{xy}) = \Lambda^{-1}$. We need to evaluate $\NDint
G\,dX$.  First separate out the strutless part, $P$, using a standard
trick: (We note that in the first line below we slightly extend the
definition of $\langle\cdot,\cdot\rangle_{\bar{X}}$, allowing it  to
have values in $\calA(X)$ rather than just $\calA(\emptyset)$ as in the
original Definition~\ref{def:NDint})
\begin{align*}
\NDint\!G\,dX &= \NDint 
  \bigg\langle P\big/(x \mapsto \prtl{\bar x})\comma
	  \exp \biggl(\frac12\,l_{xy} \strutn{x}{y} +
		       \strutv{x}{\bar x} \biggr)
  \bigg\rangle_{\!\!\bar X}\,dX \\
	&= 
  \bigg\langle P\big/(x \mapsto \prtl{\bar x})\comma
	  \NDint\!\exp \biggl(\frac12 l_{xy} \strutn{x}{y} +
			       \strutv{x}{\bar x} \biggr)\,dX
  \bigg\rangle_{\!\!\bar X}
\end{align*}
Now we complete the square in the integral:
\[
\NDint\!\exp \biggl(\frac12\,l_{xy} \strutn{x}{y} +
		\strutv{x}{\bar x} \biggr)\,dX
  =
\exp\biggl(-\frac12\, l^{xy} \strutu{\bar x}{\bar y}\biggr)
\NDint\!\exp \biggl(\frac12\, l_{xy} \strutn{x}{y} +
	\strutv{x}{\bar x} +
	\frac12\, l^{xy} \strutu{\bar x}{\bar y}\biggr)\,dX
\]
A short computation (left to the reader) shows that we have, indeed,
completed the square:
\[
\exp \biggl(\frac12\, l_{xy} \strutn{x}{y} +
	\strutv{x}{\bar x} +
	\frac12\, l^{xy} \strutu{\bar x}{\bar y}\biggr) =
  \exp\biggl(\frac12\, l_{xy} \strutn{x}{y}\biggr)
	\bigg/(x\mapsto x+l^{xy}\bar y).
\]
But now, by Propositions~\ref{prop:transinv} and~\ref{prop:gaussint},
\[
  \NDint\!\exp \biggl(\frac12\, l_{xy} \strutn{x}{y}
    + \strutv{x}{\bar x}
    + \frac12\, l^{xy} \strutu{\bar x}{\bar y}\biggr)\,dX
  = \NDint\!\exp\biggl(\frac12\, l_{xy} \strutn{x}{y}\biggr)\,dX
  = (-1)^{nm} (\det \Lambda)^m
\]
and so
\begin{align*}
  \NDint\!G\,dX &= \bigg\langle P\big/(x \rightarrow \prtl{\bar x})\comma
    \exp\biggl(-\frac12\, l^{xy} \strutu{\bar x}{\bar y}\biggr)
    \cdot (-1)^{nm} (\det \Lambda)^m \bigg\rangle_{\bar X} \\
  &= (-1)^{nm} (\det \Lambda)^m \Gint\!G\,dX \Bigg/ (O_m, C_{2m+1})
\end{align*}
\end{proof}

\begin{remark}
As in the classical case (see, e.g., the Appendix of~\cite{Aarhus:I}),
this proof can be recast in the language of Laplace (or Fourier)
transforms.
\end{remark}

\section{Proposition~\ref{prop:main} implies Theorem~\ref{thm:main}}
\lbl{sec:dull}

Let $L$ be an $n$-component regular link having linking
matrix $\Lambda$ having $\sigma_+$ positive eigenvalues and $\sigma_-$
negative eigenvalues. Lemma~\ref{lem:equivalence} and
Proposition~\ref{prop:main} imply that in the quotient 
$\Ao_{\le m}(\emptyset)/(O_m, P_{m+1})\simeq\calA_{\le m}(\emptyset)$ we
have:
\begin{equation} \lbl{eq:dull}
  \Aa_0(L) = (-1)^{nm}(\det \Lambda)^{-m} \LMO_0^{(m)}(L).
\end{equation}
Applying this equality to $U_x^\pm$, the $x$-labeled unknot with $\pm
1$ framing, we find that
\[
  \Aa_0(U_x^+) = (-1)^m \LMO_0(U_x^+)
  \qquad\text{and}\qquad
  \Aa_0(U_x^-) = \LMO_0(U_x^-).
\]
These are the renormalization factors used in the definition of $\Aa$
and $\LMO^{(m)}$, respectively. Using them and Equation~\eqref{eq:dull}
once again, we can compare $\Aa(M)$ and $\LMO^{(m)}(M)$ as follows:
\begin{align*}
  \Aa(M) &= \Aa_0(U_x^+)^{-\sigma_+} \Aa_0(U_x^-)^{-\sigma_-}
    \Aa_0(L) \\
  &= (-1)^{(n-\sigma_+)m}(\det\Lambda)^{-m}
    \LMO_0(U_x^+)^{-\sigma_+}
    \LMO_0(U_x^-)^{-\sigma_-}
    \LMO_0(L) \\
  &= |\det\Lambda|^{-m}\,\LMO^{(m)}(M) = |H^1(M)|^{-m}\,\LMO^{(m)}(M)
\end{align*}

Equation~\eqref{eq:main} of Theorem~\ref{thm:main} now follows from
the fact that $\LMO(M)$ takes its degree $m$ part from $\LMO^{(m)}(M)$.

\section{Some philosophy}\lbl{sec:philosophy}

The impatient mathematical reader may skip this section; there is nothing
with rigorous mathematical content here.  For the moment, the material in
this section is purely philosophy, and not very well-developed philosophy
at that.

This interpretation of $\NDint$ as negative-dimensional formal integration
probably seems somewhat strange.  After all, Chern-Simons theory,
the basis for the theory of trivalent graphs (the spaces $\calA$ and
$\calB$), and much of the theory of Vassiliev invariants, takes place
very definitely in positive dimensions: the vector space associated to a
vertex is some Lie algebra $\frakg$.  To integrate over these positive
dimensional spaces, a different theory is necessary, as developed in
Part II of this series~\cite{Aarhus:II}.

One potential answer to this problem is to forget about the Lie
algebra for the moment and just look at the structure of diagrams.
There are at least three different reasonably natural sign conventions
for the diagrams under consideration \cite{Kontsevich:FeynmanDiagrams,
Thurston:IntegralExpressions}.  The standard choice is to give an
orientation (ordering up to even permutations) of the edges around each
vertex.  But another natural (though usually less convenient) choice is to
leave the edges around a vertex unordered and, instead, give a direction
on each edge and a sign ordering of the {\em set} of all vertices.%
\footnote{
  In this discussion, we assume that all vertices of the graphs have
  odd valency (as holds for all diagrams considered in this paper).
  See \cite{Kontsevich:FeynmanDiagrams, Thurston:IntegralExpressions}
  for details on dealing with diagrams with vertices of even valency.
}
But now look what happens to the space $\calB$: because of the ordering
on the vertices, the diagrams are no longer completely symmetric
under the action of permuting the legs; they are now completely
{\em anti-symmetric}.  This anti-symmetry of legs is exactly what we
would expect for functions of fermionic variables or functions on a
negative-dimensional space.  Furthermore, the integration map $\NDint$
is quite suggestive from this point of view: it looks like evaluation
against a top exterior power of a symplectic form on a vector space,
which is a correct analogue of integration.

Alternatively, we could try to keep the connection
with physics, but try to find a physical theory that
exhibits this negative-dimensional behavior. Fortunately,
such a theory has been found: it is the Rozansky-Witten
theory~\cite{RozanskyWitten:HyperKahler,Kontsevich:RozanskyWitten,%
Kapranov:RozanskyWitten}.  In this theory weight systems are constructed
from a Hyper-K\"ahler manifold $Y$ of dimension $4m$.  The really
interesting thing for present purposes is that the factors assigned to the
vertices are holomorphic one-forms on $Y$, which anti-commute (using the
wedge product on forms).  (In keeping with the above remarks about signs,
each edges is assigned a symplectic form on $Y$, which is anti-symmetric.)
So in this case, there {\em is} a kind of $(-2m)$-dimensional space
associated to vertices.  (But note that this space is ``spread out''
over $Y$: it is the parity-reversed holomorphic tangent bundle.)

Finally, it's interesting to note that the definition of $\NDint$ is
more general than that of $\Gint$, and the proofs are equally simple.
On the other hand, $\Gint$ has some advantages. It is easier to
compute, as you can see in Appendix~\ref{sec:compexample}.  Its
philosophical meaning is much clearer and it takes values in
$\calA(\emptyset)$ directly, rather than in some quotient. Also, it makes
Part~IV of this series possible --- its relationship with Lie algebras
is clearer.


\appendix
\section{A computational example}
\label{sec:compexample}

In order to make more concrete the definitions of the two types of
integrals we consider, in this appendix we will integrate a formal
power series with $\NDint$ and with $\Gint$ and check that the answer
obeys Proposition~\ref{prop:main}.  We will also see how the
combinatorial factors work out in practice.

We will integrate the following ``function'' in $\calB(x,y)$:
\[
f(x,y) = \exp\left(\mathgraph{draws/comp.0} + \mathgraph{draws/comp.1}
  + \frac{1}{2} \mathgraph{draws/comp.2}\right).
\]
Since our drawings will get crowded, we will replace the labels $x$
and $\partial/\partial x$ by a solid circle~$\bullet$, and the labels $y$ and $\partial/\partial y$ by a
open circle~$\circ$.  We will write $I^{FG}$ for the formal Gaussian
integral of~$f$, and $I^{(m)}$ for the negative-dimensional integral
of~$f$.

\subsection{Computing the formal Gaussian integral}
\label{sec:comp-gint}

In order to compute the first few terms of $I^{FG} = \Gint f(x,y)\,dx\,dy$,
first note that $f$ takes the form of a perturbed Gaussian
\begin{align*}
f &= P\times\exp\biggl(l_{ij}\strutn{i}{j}\biggr) \\
  &= \exp\left(\mathgraph{draws/comp.10}\right)\times
     \exp\left(\mathgraph{draws/comp.11} + \frac{1}{2} \mathgraph{draws/comp.12}\right)
\end{align*}
with quadratic part $\Lambda = (l_{ij}) = \bigl(\begin{smallmatrix} 0&1 \\ 1&1
\end{smallmatrix}\bigr)$.  Then the prescription for the formal Gaussian
integral of $f$ involves pairing the perturbation with a quadratic part given by $-\Lambda^{-1} =
\bigl(\begin{smallmatrix} 1&-1\\-1&0 \end{smallmatrix}\bigr)$:
\[
I^{FG} =
\left\langle\exp\left(\mathgraph{draws/comp.10}\right),
     \exp\left(-\mathgraph{draws/comp.11} + \frac{1}{2} \mathgraph{draws/comp.13}\right)
   \right\rangle
\]
Expanding the exponential on the left of the pairing, we get
\[
I^{FG} =
\left\langle1,
       \exp\left(-\mathgraph{draws/comp.11} + \frac{1}{2} \mathgraph{draws/comp.13}\right)
     \right\rangle + 
     \left\langle\mathgraph{draws/comp.10},
       \exp\left(-\mathgraph{draws/comp.11} + \frac{1}{2} \mathgraph{draws/comp.13}\right)
     \right\rangle + \cdots
\]
Each left side has a definite number of vertices of each color and has a
non-zero pairing with only one term in the exponential on the right,
so we can simplify this to
\[
I^{FG} =
\langle1,1\rangle + \left\langle\mathgraph{draws/comp.10},
           -\frac{1}{2}\mathgraph{draws/comp.14}
         \right\rangle + \cdots
\]
Now we evaluate the pairing.
Note how the combinatorics of the pairing with exponentials of struts
works out: we end up summing over all ways of pairing the end points
of $P$, with coefficients given by products of appropriate entries of
$-\Lambda^{-1}$.  All other combinatorial factors cancel.
\begin{align*}
I^{FG} &=
   1
     -\mathgraph{draws/comp.20} -\mathgraph{draws/comp.21} -\mathgraph{draws/comp.22}
     + \cdots\\
  &= 1 -2\doubletheta -\crossedtet + \cdots.
\end{align*}

\subsection{Computing the negative-dimensional integral $\NDint$}
\label{sec:comp-ndint}

We now turn to the negative-dimensional integral of $f$.  To be
able to compare the results, we must have $m$ at least $2$; let us
therefore set $m=2$.
\begin{align*}
I^{(2)} &=
  \left\langle f ,
    \frac{1}{2!}\left(\frac{1}{2}\mathgraph{draws/comp.13}\right)^2\times
    \frac{1}{2!}\left(\frac{1}{2}\mathgraph{draws/comp.12}\right)^2
  \right\rangle\\
  &= \left\langle\exp\left(\mathgraph{draws/comp.10} + \mathgraph{draws/comp.11}
           + \frac{1}{2} \mathgraph{draws/comp.12}\right) ,
       \frac{1}{64}\mathgraph{draws/comp.15}\right\rangle
\end{align*}
The right side has four $x$-vertices and four $y$-vertices.  There
are only two terms in the exponential on the left which have non-zero
pairing with the right side:
\[
I^{(2)} =
  \left\langle\frac{1}{24}\mathgraph{draws/comp.50} + \frac{1}{2}\mathgraph{draws/comp.60},
       \frac{1}{64}\mathgraph{draws/comp.15}\right\rangle
\]
To evaluate the pairing, notice that the
combinatorial factors on the right side cancel out, and we end up
summing over all ways of pairing the vertices of the same color on the
left, each appearing with coefficient 1.
\begin{align*}
I^{(2)} &=
   \frac{1}{24}\left(
      \mathgraphsmall{draws/comp.51}+\mathgraphsmall{draws/comp.52}+\mathgraphsmall{draws/comp.53}
     +\mathgraphsmall{draws/comp.54}+\mathgraphsmall{draws/comp.55}+\mathgraphsmall{draws/comp.56}
     +\mathgraphsmall{draws/comp.57}+\mathgraphsmall{draws/comp.58}+\mathgraphsmall{draws/comp.59}
     \right) \\
   &\qquad+
   \frac{1}{2}\left(
      \mathgraphsmall{draws/comp.61}+\mathgraphsmall{draws/comp.62}+\mathgraphsmall{draws/comp.63}
     +\mathgraphsmall{draws/comp.64}+\mathgraphsmall{draws/comp.65}+\mathgraphsmall{draws/comp.66}
     +\mathgraphsmall{draws/comp.67}+\mathgraphsmall{draws/comp.68}+\mathgraphsmall{draws/comp.69}
     \right)
\end{align*}
Now apply the $O_m$ relation with $m=2$, which replaces each circle by
$-4$.
\begin{align*}
I^{(2)} &=
  \frac{1}{24} \Bigl(16-4-4-4+16-4-4-4+16\Bigr)\\
   &\qquad+\frac{1}{2}\Bigl(-4\doubletheta+\doubletheta+\doubletheta
                     -4\crossedtet+\crossedtet+\crossedtet
                     -4\doubletheta+\doubletheta+\doubletheta\Bigr) \\
  &= 1 - 2\doubletheta - \crossedtet.
\end{align*}

Note that this matches the result we found from the formal Gaussian
integration, which is expected since $\det \Lambda = -1$ and the degree of
the non-trivial elements is even.


\section*{Acknowledgement}

The seeds leading to this work were planted when the four of us (as
well as Le, Murakami (H\&J), Ohtsuki, and many other like-minded
people) were visiting \Arhus, Denmark, for a special semester on
geometry and physics, in August 1995. We wish to thank the organizers,
J.~Dupont, H.~Pedersen, A.~Swann and especially J.~Andersen for their
hospitality and for the stimulating atmosphere they created.  We also
wish to thank N.~Habegger, M.~Hutchings, T.~Q.~T.~Le J.~Lieberum,
A.~Referee and N.~Reshetikhin for additional remarks and suggestions,
the Center for Discrete Mathematics and Theoretical Computer Science at
the Hebrew University for financial support, and the
Volkswagen-Stiftung (RiP-program in Oberwolfach) for their hospitality
and financial support.

\ifx\undefined\bysame
        \newcommand{\bysame}{\leavevmode\hbox to3em{\hrulefill}\,}
\fi

\end{document}